\documentclass[letterpaper, 9pt, conference]{ieeeconf}  

\IEEEoverridecommandlockouts{}                               
\usepackage{mw88-latex-cdc-customizations}
\usepackage{tikzscale}
\usepackage{proof-at-the-end}

\pgfkeys{/prAtEnd/global custom defaults/.style={
        text link={\noindent
                The \hyperref[proof:prAtEnd\pratendcountercurrent]{proof} can be found in the Appendix on page~\pageref{proof:prAtEnd\pratendcountercurrent}.},
        text proof={ of
                \string\Cref{thm:prAtEnd\pratendcountercurrent}},
        }
}

\usepackage[utf8x]{inputenc}

\makeatletter
\declaretheoremstyle[
        headpunct={},
        notebraces={\lparen}{\rparen},
        spaceabove=1em,
        spacebelow=1em,
        headformat={\NOTE},
        qed=\(\lrcorner\),
    ]{customstyle}
\declaretheorem[
        name=cass,
        numbered=no,
        style=customstyle,
        Refname={Assumption,Assumptions},
    ]{customassumption}
\makeatother

\title{\LARGE \bf
Inferring Global Exponential  Stability Properties using Lie-bracket Approximations (Extended Version)
}

\author{Marc Weber, Bahman Gharesifard and Christian Ebenbauer
\thanks{M. Weber and C. Ebenbauer are with Faculty of Electrical Engineering and Information Technology,
        RWTH Aachen University, D-52074 Aachen, Germany
        {\tt\small weber/ebenbauer@ic.rwth-aachen.de}}%
\thanks{B. Gharesifard is with the Electrical \& Computer Engineering Department, University of California, Los Angeles, USA
        {\tt\small gharesifard@ucla.edu}}%
}

\begin{document}

\maketitle
\thispagestyle{empty}
\pagestyle{empty}

\begin{abstract}
  In the present paper, a novel result for inferring uniform global, not semi-global, exponential stability in the sense of Lyapunov with respect to input-affine systems from global uniform  exponential stability properties with respect to their associated Lie-bracket systems is shown.
  The result is applied to adapt dither frequencies to find a sufficiently high gain in adaptive control of linear unknown systems, and a simple numerical example is simulated to support the theoretical findings.
\end{abstract}




\section{Introduction}\label{sec:problem-statement}
 Consider the input-affine system
\begin{align}
  \dot{x}(t) & = f_0(t, x(t)) + \sum_{i=1}^l \omega^{p_i} f_i(t, x(t)) \, u_i(\omega t),
  \tag{S}
  \label{eqn:es-system}
\end{align}
with \( x(t) \in \mathbb{R}^n \) the state vector of the system, \( u_i(t) \in \mathbb{R} \) the control inputs, \( \omega \in \interval[open]{0}{\infty}\) and \( p_i \in \interval[open]{0}{1} \).
Under certain assumptions, e.g.~\cite[Assumptions~1.1-1.3]{labar2019thesis}, a so-called \emph{Lie-bracket system}~\eqref{eqn:lbs} may be associated with~\eqref{eqn:es-system}, namely
\begin{subequations}
  \begin{align}
    \dot{\bar{x}}(t) & = f_0(t, \bar{x}(t)) +
    \smash{
        \lim_{\omega \to \infty} \sum_{\begin{limitarray}
                                             1 &\leq i &< l \\
                                             i &< j &\leq l
                                           \end{limitarray}}
      } \commutator{f_i}{f_j}(t, \bar{x}(t)) \, \gamma_{ij}(\omega),
    \tag{LBS}
    \label{eqn:lbs}
  \end{align}
  with \( T \coloneq 2\pi \omega^{-1}\) and
  \begin{align}
    \gamma_{ij}(\omega) & \coloneq \frac{\omega^{p_i + p_j}}{T} \int_{0}^T \int_0^\theta u_j(\omega \theta) \, u_i(\omega \tau) \odif{\tau} \odif{\theta}.
    \label{eqn:gamma-ij}
  \end{align}
\end{subequations}

Depending on the choice of input vector fields~\( f_1, \ldots, f_l\), the Lie-bracket system potentially exhibits additional control directions not contained in the distribution of \( f_1, \ldots, f_l\), which may be exploited to steer the state in a desired direction.

Using Lie-bracket systems, it has formally been shown in~\cite{DURR20131538} and~\cite{labar2019thesis}, that for finite intervals in time, solutions \( x_{\omega}(t; t_0, x_0)\) to~\eqref{eqn:es-system} with initial time \( t_0\) and initial condition~\(x_0\), parametrized by \(\omega\),  can be made to approximate solutions~\( \bar{x}(t; t_0, x_0)\) with identical initial time and initial conditions up to an arbitrary but constant user-specified error bound, by choosing~\( \omega\) sufficiently large.

These properties are then used to obtain stability results, where from (global) uniform asymptotic stability of the origin of~\eqref{eqn:lbs},  (semi-global) practical uniform asymptotic stability for the origin of~\eqref{eqn:es-system} is inferred.
We refer to the definitions in~\cite[Definition 1-5]{DURR20131538} for the concept of practical stability.

The objective of the present work is to provide tools for transferring {global} stability properties in the sense of Lyapunov from associated Lie-bracket systems over as \textbf{global} stability properties to their input-affine counterparts for a suitably restricted class of input vector fields.

Systems of form~\eqref{eqn:es-system} arise most prominently in the field of \emph{Extremum Seeking}, where they are usually coupled with some output function, with the goal of steering the state~\(x\) in a model-free way towards a minimizer of the output in real-time without a-priori knowledge of the minimizer.
A comprehensive history of extremum seeking is given in~\cite{scheinker2017model}, where its origins from adaptive control and fundamental proofs by singular perturbations are described.
Since \num{2010}, there has been a trend to employ extremum seeking for stabilization~\cite{scheinker2017model}.
With a shift of focus from optimizing the output, to driving a control Lyapunov function to zero, the singular perturbations approach used throughout historic proofs has found its contemporary counterparts in Lie-bracket approximation interpretation~\cite{DURR20131538}\cite{DURR20142591}\cite{SUTTNER2019214}.
The most prominent results, however, lead at most to semi-global practical uniform asymptotic stability~\cite{DURR20131538}\cite{Drr2017ExtremumSF}.
In~\cite{Grushkovskaya2019PartialStabilityConcept}, a Lyapunov-like theorem is proven for inferring the practical stability properties of the origin for~\eqref{eqn:es-system}.
Moreover, the authors provide sufficient conditions and examples for pairs of vector fields generating a gradient-descent-like Lie-bracket system.
Lastly, they also provide an application example for practically stabilizing a given control system using a Control-Lyapunov Function and a high-gain controller.
One key point of the results in \cite{Grushkovskaya2019PartialStabilityConcept} is, that the notion of practical stability is extended in a way that enables statements about boundedness of additional co-states.
{Global} exponential stability in the sense of Lyapunov has been shown in~\cite{SUTTNER201715464}, however through the use of  adaptive dither frequencies.
A recent result after the submission of this paper in \cite[Theorem~2(b)]{abdelgalil2024initializationfreeliebracketextremumseeking} provides a global uniform practical asymptotic stability result with slightly milder assumptions than our current approach.
To the best of our knowledge, no other work has targeted {global} stability properties.
In our previous work~\cite{weber2021note}, we investigated the possibility of adaptive high-gain stabilization for completely unknown linear systems through Lie-bracket approximation, with the goal of combining the vibrational stabilization from~\cite{Grushkovskaya2018GeneratingVectorFields} with a self-tuning gain.
Since the state may diverge greatly and to an a-priori unknown value from the origin during the self-tuning, we additionally require more than just semi-global stability properties.

The current work presents our results on inferring \textbf{global} stability properties in the sense of Lyapunov for~\eqref{eqn:es-system} from stability properties of its associated Lie-bracket system in contrast to the common semi-global stability results in literature.
This stronger result comes at the cost of stronger assumptions, namely global Lipschitz properties of the vector fields and Lie derivatives involved.
In contrast to~\cite{Grushkovskaya2019ExponentialStabilization}, our analysis does not use Lyapunov-function arguments, but is directly based on trajectory approximation properties in the spirit of~\cite{labar2019thesis}.

\section{Main Results}\label{sec:main-results}
 \noindent
Lie-bracket approximation utilizes the non-commutativity of flows of the input vector fields~\(f_i\) to generate additional control directions, which could otherwise not be obtained by linearly combining the existing input vector fields. To do so, \Cref{ass:ui} will introduce the so-called \emph{virtual inputs}~\(u_i\) used to excite the control directions~\(f_i\), by periodically switching between their respective flows. The Lie-bracket as a commutator of vector fields is then exploited to provide the desired new control direction.
\Cref{ass:fi,ass:pi} are used on the remainder expressions arising during functional expansion of solutions to show that the undesired remainder terms in our calculations will tend to zero.

In the sequel, we require the following assumptions:
\begin{assumption}[Conditions on \( u_i \)]\label{ass:ui}
  For all \( i \in \set{ 1, \dots, l} \),
  \begin{enumerate}[
      ref={\theassumption.\arabic*},
    ]
    \item\label[assumption]{ass:ui:1} \( u_i: \mathbb{R} \to \mathbb{R} \) is a measurable function, such that \( \sup_{t \in \mathbb{R}} |u_i(t)| \leq 1 \);
    \item\label[assumption]{ass:ui:2} \( u_i \) is \(2\pi\)-periodic, i.e.\ \( u_i(t + 2\pi) = u_i(t)\) for all \( t \in \mathbb{R} \);
    \item\label[assumption]{ass:ui:3} \( u_i \) has zero mean on one period, i.e. \!\( {\int_{0}^{2\pi} \!\!\!\!\! u_i(t) \odif{t} = 0} \). \qedhere
  \end{enumerate}
\end{assumption}
The \emph{virtual inputs}~\( u_i \), also called \emph{dithers} or \emph{dither signals} are bounded by an arbitrary constant, which we took to be \num{1} without loss of generality.
Observe, that \Cref{ass:ui:3} prevents~\eqref{eqn:es-system} from being represented as a driftless system by moving \(f_0\) into the summation with \(u_0(\omega t) = 1\).
Periodicity and zero mean ensure that all different control directions are explored equally and no additional drift except the vector field~\( f_0 \) enters the system.
\begin{assumption}[Conditions on \( f_i \)]\label{ass:fi}
  \noindent
  There exists an \( L \in \interval[open]{0}{\infty} \), such that for all \({ i,j,m\in \set{0, 1, \ldots, l} } \),
  {
      \setlength{\multicolsep}{\abovedisplayskip}
      \begin{multicols}{3}
        \begin{enumerate}[
            ref={\theassumption.\footnotesize{\arabic*}},
          ]
          \item\label[assumption]{ass:fi:2:1} \( \phantom{\pdif{1}} f_i \),
          \item\label[assumption]{ass:fi:2:4} \( \pdif{1} f_i \),
          \item\label[assumption]{ass:fi:2:2} \( \phantom{\pdif{1}} \mathcal{L}_{f_j}f_i\),
          \item\label[assumption]{ass:fi:2:5} \( \pdif{1} \mathcal{L}_{f_j}f_i\),
          \item\label[assumption]{ass:fi:2:3} \( \mathcal{L}_{f_m}\mathcal{L}_{f_j}f_i  \)
        \end{enumerate}
      \end{multicols}
    }
  \noindent
  are continuous in the first argument, {globally uniformly \(L\)-Lipschitz} in the second argument and defined for all \( t \in \mathbb{R} \), and almost every \(x \in \mathbb{R}^n\),
  Additionally, we require
  \begin{align*}
    f_i(t, 0)                            & = 0, & \mathcal{L}_{f_j}f_i (t, 0) & = 0,      & \mathcal{L}_{f_m}\mathcal{L}_{f_j}f_i(t, 0) & = \mathrlap{0.} \\
    \pdif{1} f_i(t, 0)                   & = 0, & 
    \pdif{1} \mathcal{L}_{f_j}f_i (t, 0) & = 0, &                             & \qedhere
  \end{align*}
\end{assumption}
Since we will rely on functional expansion for our arguments, the \emph{control vector fields} need to be sufficiently differentiable.
\Cref{ass:fi:2:1,ass:fi:2:3} imply that the vector fields~\( f_i \) and \(\mathcal{L}_{f_j}{f_i}\) are absolutely continuous, and hence almost everywhere differentiable.
The global Lipschitz conditions are technical and restrictive; however they are needed throughout our arguments to arrive at \textbf{global} stability properties.
Nevertheless, we show by constructing an example that the resulting class of vector field tuples is not empty.

\begin{assumption}[Interaction condition]\label{ass:pi}
  \\
  \noindent
  For all \( {i, j \in \set{ 1,2, \ldots, l }} \), the vector fields~\(f_i\) and virtual inputs~\(u_i\) satisfy the interaction condition:\\
  If \( p_i + p_j > 1\), then either
  \begin{enumerate}[
      ref={\theassumption.\arabic*},
    ]
    \item \( \int_{0}^{T} \int_{0}^{\theta} u_j(\omega \theta) \, u_i(\omega \tau) \odif{\tau} \odif{\theta} = 0\), or
    \item \(\commutator{f_i}{ f_j}(t, x) = 0\) for all \( x \in \mathbb{R}^n\) and \( t \in \mathbb{R}\).
          \qedhere
  \end{enumerate}
\end{assumption}
Finally, \Cref{ass:pi} describes requirements on the interaction between the \emph{powers}~\( p_i \) \emph{of the dither frequencies} and the Lie-brackets of control vector fields. The assumptions on \(p_i, i \in \set{1, \ldots, l}\) are needed for the proofs;
We take in practice and in our simulations \( p_i=\tfrac 1 2\) for all \(i \in \set{1, \ldots, l}\), such that \Cref{ass:pi} does not impose any constraints.

With these assumptions established, we proceed to state our main results.

Firstly, we give a modified version of the trajectory approximation lemma in~\cite[Lemma~1.1]{labar2019thesis}, relating the approximation error to the norm of the initial condition instead of providing absolute bounds.
\begin{theorem}[Trajectory approximation on finite intervals]\label{lemma:approximation-lemma}
  Consider a system of form~\eqref{eqn:es-system} with \Cref{ass:ui,ass:fi,ass:pi} satisfied.
  Let
  \( t_f \in \interval[open]{0}{\infty}\)
  and \( \bar\Lambda \in \interval[open right]{1}{\infty}\), such that,
  for every \( t_0 \in \mathbb{ R}\) and every \( x_0 \in \mathbb{R}^n\),
  the trajectory of~\eqref{eqn:lbs} through \( \bar{x}(t_0) = x_0 \) exists, is unique and satisfies
  \begin{align}
    \norm{ \bar{x}_{\omega}(t; t_0, {x}_0) } & \leq \bar\Lambda \, \norm{{x}_0 }\label{eqn:lemma:approximation-lemma:condition}
  \end{align}
  for all \( t \in \interval[]{t_0}{t_0+t_f}\).
  
  Then, for every \( D \in  \interval[open]{0}{\infty}\), there exists an \( \omega^\star \in \interval[open]{0}{\infty}\), such that, for all \( \omega \in \interval[open]{\omega^\star}{\infty}\), for all \( t_0 \in \mathbb{R}\) and for all \( x_0 \in \mathbb{R}^n \!\setminus \!\set*{0}\), the trajectories of~\eqref{eqn:lbs} and~\eqref{eqn:es-system} through \( x(t_0) = \bar{x}(t_0) = x_0 \) satisfy
  \begin{align}
    \norm{ x_\omega(t; t_0, x_0) - \bar{x}_\omega(t; t_0, x_0)} & < D \, \norm{x_0}    \label{eqn:lemma:approximation-lemma:approximation}
  \end{align}
  for all \( t \in \interval[]{t_0}{t_0+t_f}\). In particular, define
  \begin{subequations}
    \begin{align*}
      \mathcal{P}_1
                    & \coloneq \set*{          (i, j) \in \set{1, \ldots, l}^2 \given p_i + p_j < 1}         \\
      \mathcal{P}_2 & \coloneq \set*{          (i, j, m) \in \set{1, \ldots, l}^3 \given p_i + p_j +p_m < 2}
    \end{align*}
    and
    \begin{align}
      p'   & \coloneq \quad\smashoperator{\max_{i\in \set{1, \ldots, l}}} \quad p_i, \label{eqn:def:p:1}              \\
      p''  & \coloneq \quad\smashoperator{\max_{(i,j) \in \mathcal{P}_1}}\quad p_i + p_j, \label{eqn:def:p:2}         \\
      p''' & \coloneq \quad\smashoperator{\max_{ (i,j, m) \in \mathcal{P}_2}}\quad p_i + p_j + p_m\label{eqn:def:p:3}
    \end{align}
    and let
    \begin{align}
      p^\star & \coloneq \min\set*{ 1- p',  1-p'',  2- p'''}.
      \label{eqn:def:p-star}
    \end{align}
    Then,
    \begin{align}
      \omega^\star & \coloneqq
      \max\set*{ 1,
          \sqrt[p^\star]{ c \tfrac{\bar\Lambda}{D} \left\lparen
              9t_f +12 +8 \pi
              \right\rparen e^{2 \pi l^2 L t_f}
            }
        }
      \label{eqn:def:omega-star}
    \end{align}
    with \( c \coloneq \pi^2 {(l+1)}^3 L^2\) is a sufficient choice. \qedhere
  \end{subequations}
\end{theorem}

The \hyperref[proof:approximation-lemma]{proof} can be found in the Appendix on page~\pageref{proof:approximation-lemma}.

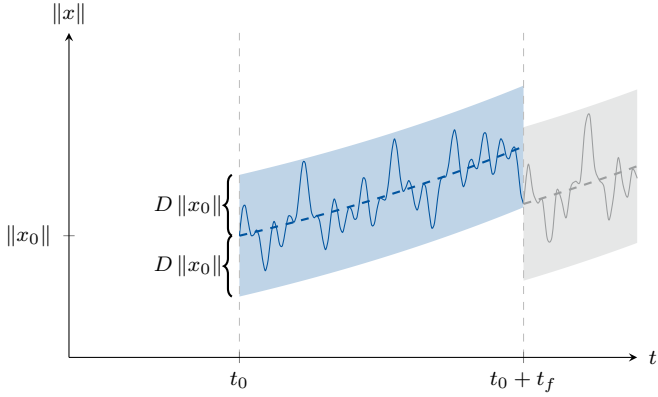
\begin{figure}
  \centering%
  \pgfplotsset{trig format plots=rad, compat=1.12}
\tikzmath{
    real \initialCondition;
    real \initialTime;
    real \tf;
    real \D;
    \initialCondition = 0.3;        
    \initialTime = 0.3;      
    \tf = 0.5;       
    \D = 0.5;
  }
\def\lbs(#1, #2, #3){
    (#2)*exp(1.1*((#1)-(#3)))
  }
\def\ai(#1, #2, #3){
    (#2)*exp(1.1*((#1)-(#3)))
    +
    \D*(#2)*sin(97*((#1)-(#3)))*cos(57*((#1)-(#3)))*cos(82*((#1)-(#3)))
  }
\tikzmath{
    real \specialY;
    \specialY = 3.7835699e-01;
  }
%
\makeatletter
\let\real=\pgfmath@calc@real
\let\minof=\pgfmath@calc@minof
\let\maxof=\pgfmath@calc@maxof
\let\ratio=\pgfmath@calc@ratio
\let\widthof=\pgfmath@calc@widthof
\let\heightof=\pgfmath@calc@heightof
\let\depthof=\pgfmath@calc@depthof
\makeatother
\begin{tikzpicture}
  \begin{axis}[
      axis on top,
      scale only axis,
      axis x line = center,
      axis y line = center,
      xmin = 0,
      xmax = 1,
      ymin = 0,
      ymax = 0.8,
      width = {\linewidth-\widthof{\(\norm{x_0}\)}-0.3cm-\widthof{t}},
      height = 0.5\linewidth,
      xlabel style = {
          inner sep = 0pt,
          align = right,
          xshift=0.15cm
        },
      yticklabel style = {
          inner sep = 0pt,
          align = right,
          xshift=-0.15cm
        },
      xlabel = {\(t\)},
      ylabel = {\(\norm{x}\)},
      extra x ticks={\initialTime, (\initialTime+\tf)},
      extra x tick labels={\(t_0\), \(t_0+t_f\)},
      extra y ticks = {\initialCondition},
      extra y tick labels = {\(\norm{x_0}\)},
      extra x tick style = { grid = major, },
      ytick=\empty,
      xtick=\empty,
      clip=false,
    ]
    \addplot[
        name path=UpperBound,
        color=black,
        domain=\initialTime:\initialTime+\tf,
        samples=20,
        smooth,
        draw=none,
      ] {\D*\initialCondition + \lbs(x, \initialCondition, \initialTime) }; \label{plot:trajectory-approximation:upper-bound}
    \addplot[
        name path=LowerBound,
        color=black,
        domain=\initialTime:(\initialTime+\tf),
        samples=20,
        smooth,
        draw=none,
      ] {-\D*\initialCondition + \lbs(x, \initialCondition, \initialTime)}; \label{plot:trajectory-approximation:lower-bound}
    \addplot[blue!25] fill between[of=LowerBound and UpperBound];\label{plot:trajectory-approximation:error-interval}
    \addplot[
        color=blue,
        dashed,
        domain=\initialTime:(\initialTime+\tf),
        samples=20,
        smooth,
      ] { \lbs(x, \initialCondition, \initialTime) }; \label{plot:trajectory-approximation:lbs}
    \addplot[
        color=blue,
        thin,
        domain=\initialTime:(\initialTime+\tf),
        samples=100,
        smooth,
      ] { \ai(x, \initialCondition, \initialTime) }; \label{plot:trajectory-approximation:ai}
    
    \addplot[
        color=black50,
        dashed,
        domain=(\initialTime+\tf):1,
        smooth,
      ] { \lbs(x, \specialY, \initialTime+\tf)}; \label{plot:trajectory-approximation:lbs-continuedc}
    \addplot[
        color=black50,
        thin,
        domain=(\initialTime+\tf):1,
        smooth,
      ] { \ai(x, \specialY, \initialTime+\tf)}; \label{plot:trajectory-approximation:ai-continuedc}
    \addplot[
        name path=UpperBoundC,
        color=black,
        domain=(\initialTime+\tf):1,
        samples=20,
        smooth,
        draw=none,
      ] {\D*\specialY + \lbs(x, \specialY, \initialTime+\tf) }; \label{plot:trajectory-approximation:upper-boundc}
    \addplot[
        name path=LowerBoundC,
        color=black,
        domain=(\initialTime+\tf):1,
        samples=20,
        smooth,
        draw=none,
      ] {-\D*\specialY + \lbs(x, \specialY, \initialTime+\tf)}; \label{plot:trajectory-approximation:lower-boundc}
    \addplot[black50!25] fill between[of=LowerBoundC and UpperBoundC];\label{plot:trajectory-approximation:error-intervalc}
    \draw [decorate, thick,  decoration={brace,amplitude=3pt}, xshift=-1mm] (axis cs:\initialTime,\initialCondition) -- (axis cs:\initialTime,{\initialCondition+\D*\initialCondition}) node [midway, anchor=east, font=\footnotesize]{\( D\norm{x_0}\)};
    \draw [decorate, thick,  decoration={brace,amplitude=3pt, mirror}, xshift=-1mm] (axis cs:\initialTime,\initialCondition) -- (axis cs:\initialTime,{\initialCondition-\D*\initialCondition}) node [midway, anchor=east, font=\footnotesize]{\( D\norm{x_0}\)};
  \end{axis}%
\end{tikzpicture}%
%
  \caption{Solutions~\(x_\omega(t; t_0, x_0)\)~(\ref{plot:trajectory-approximation:ai}) of~\eqref{eqn:es-system}, solution~\(\bar{x}(t; t_0, x_0)\)~(\ref{plot:trajectory-approximation:lbs}) of the associated Lie-bracket system~\eqref{eqn:lbs}, error bound interval~\(D\,\norm{x_0}\)~(\ref{plot:trajectory-approximation:error-interval}) over approximation horizon~\(\interval[]{t_0}{{t_0+t_f}}\) and insinuated restart at \({t_0+t_f}\) and continuation for solutions~\(x_\omega(t; {t_0+t_f}, x_\omega(t_0+t_f; t_0, x_0)), \bar{x}(t; {t_0+t_f}, x_\omega({t_0+t_f}; t_0, x_0))\) and their respective error bound interval~\(D\,\norm{x(t_f; t_0, x_0)}\)~(\ref{plot:trajectory-approximation:error-intervalc}).}\label{fig:trajectory-approximation-principle}%
\end{figure}

\Cref{lemma:approximation-lemma} states an approximation property between solutions \(x_{\omega}(t; t_0, x_0)\) and \(\bar{x}(t; t_0, x_0)\) over a finite time horizon \( \interval[]{t_0}{{t_0+ t_f}}\) of length \( t_f\).
For an illustration see \Cref{fig:trajectory-approximation-principle}.
Notably, the error bound between solutions in~\eqref{eqn:lemma:approximation-lemma:approximation} is not stated as an absolute value, but as a dependency on the norm of their common initial condition.
We can use this result to iteratively extend solutions by a length of \(t_f\) to infer \textbf{global} exponential stability properties of a system of form~\eqref{eqn:es-system} from the properties of its associated Lie-bracket system~\eqref{eqn:lbs}:

\begin{theorem}[Global uniform exponential stability]
  \label{lemma:global-exponential-stability-result}
  \mbox{}
  Consider systems of form~\eqref{eqn:es-system} with the associated Lie-bracket system~\eqref{eqn:lbs} satisfying \Cref{ass:ui,ass:fi,ass:pi}.
  Suppose that the origin is globally uniformly exponentially stable for the Lie-bracket system~\eqref{eqn:lbs}, i.e.\ there exist constants \( \bar\alpha \in \interval[open right]{1}{\infty} \) and \( \bar\beta \in \interval[open]{0}{\infty} \), such that for all \( t_0 \in \mathbb{R} \) and all \( x_0\in \mathbb{R}^n \),
  \begin{equation*}
    \norm{\bar{x} (t; t_0, x_0)}  \leq \bar{\alpha} \norm{x_0} e^{-\bar{\beta} (t -t_0)}
    \qquad\text{for all \( t \in \interval[open right]{t_0}{\infty} \).}
  \end{equation*}
  Then,
  there exists an \( \omega^\star \in \interval[open]{0}{\infty}\), such that
  the origin is  globally uniformly exponentially  stable with respect to~\eqref{eqn:es-system} for all \( \omega \in \interval[open]{\omega^\star}{\infty}\),  i.e.\
  there exist constants \( \omega^\star \in \interval[open]{0}{\infty}\), \({ \alpha \in \interval[open right]{1}{\infty}} \), \( \beta \in \interval[open]{0}{\infty} \), such that, for all \( \omega \in \interval[open]{\omega^\star}{\infty}\), for all~\( t_0\in \mathbb{R} \) and for all~\( x_0 \in \mathbb{R}^n \),
  \begin{align*}
    \norm{x_\omega(t; t_0, x_0)} & \leq \alpha \, \norm{x_0} \,  e^{-\beta (t-t_0)}
  \end{align*}
  for all \( t \in \interval[open right]{t_0}{\infty} \).
\end{theorem}

The \hyperref[proof:global-exponential-stability-result]{proof} can be found in the Appendix on page~\pageref{proof:global-exponential-stability-result}.

Observe, that~\eqref{eqn:fb1:choice-of-tf-and-d},~\eqref{eqn:lemma:ges:alpha},~\eqref{eqn:lemma:ges:beta} is one particular way to relate the global uniform exponential stability of the origin of~\eqref{eqn:es-system}  to the parameters~\(\bar\alpha\) and \(\bar\beta \) describing the global uniform exponential stability of the origin of its associated Lie-bracket system~\eqref{eqn:lbs} by choosing an error scale~\(D\) and a time horizon~\(t_f\), over which we apply the virtual inputs.
Keep in mind that a sufficiently high lower bound~\(\omega^\star\) on the dither frequency~\(\omega\) still may be determined according to~\eqref{eqn:def:omega-star}, but this same bound applies at all times and all states.
This allowed us to show genuine {global} (exponential) uniform stability in contrast to the common  semi-global practical uniform stability results outlined in the introduction.


Lastly, we present an application of the previous result.

\begin{lemma}[Global convergence with dither frequency adaptation]\label{lemma:frequency-adaptation}
  \begin{subequations}\label{eqn:es-system-frequency}
    Consider the system parametrised in \(k \in \mathbb{N}_0\)
    \begin{align}
      \dot{x}(t) & =
      f_0(t, x(t)) + \sum_{i=1}^l w_k^{p_i} f_i(t, x(t)) u_i(w_k t),
      \label{eqn:augmented-system:state}
      \\
      \dot{w}(t) & = \norm{x(t)}^2,
      \label{eqn:augmented-system:frequency}
    \end{align}
    with \( p_i = \tfrac 1 2\) for all \( i \in \set{1, \ldots, l}\) where
    \begin{align}
      w_k & \coloneq w(t_0+ kt_f; t_0, w_0)
    \end{align}
    is the solution to~\eqref{eqn:augmented-system:frequency} with initial time~\( t_0\) and initial condition~\( w_0\) sampled at a constant rate~\( t_f\).
    Suppose that the origin is globally uniformly exponentially stable for  the Lie-bracket system associated with~\eqref{eqn:augmented-system:state} for all \(k \in \mathbb{N}_0\), i.e.\ there exist constants \( \bar\alpha \in \interval[open right]{1}{\infty} \) and \( \bar\beta \in \interval[open]{0}{\infty} \) such that for all \( t_0 \in \mathbb{R} \) and all \( x_0\in \mathbb{R}^n \),
    \begin{align}
      \norm{\bar{x}(t)} & \leq \bar\alpha  \norm{x_0}  e^{-\bar\beta (t -t_0)}
      \label{eqn:fb1:adaptive-high-gain-lbs}
    \end{align}
    for all \( t \in \interval[open right]{t_0}{\infty} \). Suppose further, that \( t_f \in \interval[open]{0}{\infty}\) has been chosen sufficiently large, such that there exists a \( D \in \interval[open]{0}{\infty} \) satisfying~\eqref{eqn:fb1:choice-of-tf-and-d}.
    
    Then, for any \( x_0 \in \mathbb{R}^n\) and any \(\omega_0 \in \mathbb{R}\), \( \lim_{ t \to \infty} \norm{x(t; t_0, x_0)} = 0\), \(w(t; t_0, w_0)\) is bounded and converges to a finite limit.
  \end{subequations}
\end{lemma}

The \hyperref[proof:frequency-adaptation]{proof} can be found in the Appendix on page~\pageref{proof:frequency-adaptation}.

In \Cref{lemma:frequency-adaptation}, we showed that the problem of finding \emph{sufficiently high gains} required for the application of~\Cref{lemma:global-exponential-stability-result} can be practically adressed by adding a scalar quadratic integrator, quite similar to the Willems-controller in~\cite{willems1984global}, and sampling it, which introduces a strictly monotonically increasing sequence of dither frequencies~\(w_k\).
Taking numerical obstacles not into consideration, the proposed algorithm in \Cref{lemma:frequency-adaptation} provides a novel adaptive way of selecting a sufficient dither frequency, which may in practice lie significantly lower than the conservative bound stated in~\eqref{eqn:def:omega-star}.

\section{Simulation Example}\label{sec:simulation-example}
 We demonstrate the importance of \Cref{lemma:global-exponential-stability-result} using an example from adaptive control.

\begin{figure}
  \centering
  \begin{tikzpicture}
  \pgfplotstableread[
      skip first n = 9, 
      col sep =comma,
      header = true,
    ]{matlab/exponential_convergence.csv}{\tableAi}
  \begin{axis}[
      axis x line = center,
      axis y line = center,
      xmin = 0,
      xmax = 0.99,
      ymin= 0,
      ymax = 1.45,
      width=\linewidth-0.15cm,
      height = 0.6\linewidth,
      xlabel style = {
          inner sep = 0pt,
          align = right,
          xshift=0.15cm
        },
      xlabel = {\(t/\si{\second}\)},
      ylabel={\(x\)},
      extra y ticks={0, 1},
      extra y tick labels={0, \(x_0\)},
      ytick=\empty,
      clip=false,
    ]
    \addplot[
        color=black,
        domain=0:1,
        dotted,
        thin,
      ] {(1.39475)*exp(-1.80132*x)}; \label{plot:exponential-convergence:upperbound-analytical}
    \addplot[
        color=black,
        domain=0:1,
        dashed,
        thin,
      ] {(1.39475)*exp(-4*x)}; \label{plot:exponential-convergence:upperbound-emprical}
    \addplot[
        color=black,
        domain=0:1,
        dashed,
        thin,
      ] {(2-1.39475)*exp(-4*x)}; \label{plot:exponential-convergence:lowerbound}
    \addplot[
        thin,
        color=blue
      ] table[
        x=t,
        y = ai.x,
      ] {\tableAi};\label{plot:exponential-convergence:ai}
    \addplot[
        color=blue,
        dashed,
      ] table[
        each nth point={10},
        x=t,
        y = lbs.x,
      ] {\tableAi};\label{plot:exponential-convergence:lbs}

    \draw [
        decorate,
        thick,
        decoration={brace,amplitude=3pt},
        xshift=-1mm
      ] (axis cs:0,1) --
      node [
        midway,
        anchor=east,
        font=\footnotesize,
        inner sep=0pt,
        align=right,
        xshift=-0.15cm
      ]{\( D\norm{x_0}\)}
      (axis cs:0,1.39475)    ;
    \draw [
        decorate,
        thick,
        decoration={brace,amplitude=3pt, mirror},
        xshift=-1mm,
      ] (axis cs:0,1) --
      node [
        midway,
        anchor=east,
        font=\footnotesize,
        inner sep=0pt,
        align=right,
        xshift=-0.15cm
      ]{\( D\norm{x_0}\)}
      (axis cs:0,0.6053)
    ;
  \end{axis}
  
\end{tikzpicture}
  \caption{Simulation results show the solution~\(x_{\omega}(t; t_0, x_0)\)~(\ref{plot:exponential-convergence:ai}) from~\eqref{eqn:example-system}, the solution~\(\bar{x}(t; t_0, x_0)\)~(\ref{plot:exponential-convergence:lbs}) of the associated Lie-bracket system~\eqref{eqn:example:lbs} for a common initial condition \(x_0\), the conservative analytical upper bound~\(\alpha e^{-\beta t}\)~(\ref{plot:exponential-convergence:upperbound-analytical}) as well as the empirical  bounds~\( \alpha e^{-\bar{\beta} t}\)~(\ref{plot:exponential-convergence:upperbound-emprical}) for parameters from \Cref{sec:simulation-example}.}
  \label{fig:simulation-result}
\end{figure}
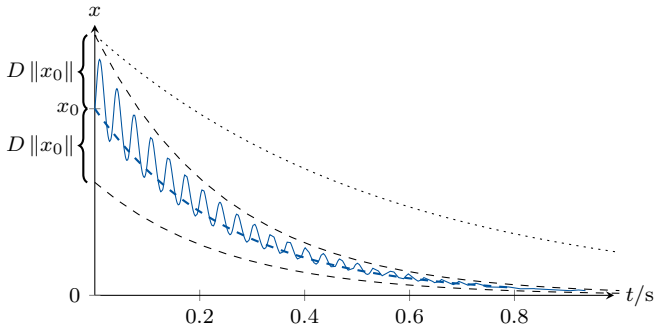

\begin{description}[
    labelindent=0pt,
    leftmargin=0pt,
    topsep=0.5em,
    itemsep=0.5em,
  ]
  \item[System Dynamics]
        Consider the scalar linear time-invariant system
        \begin{subequations}\label{eqn:example-system}
          \begin{align}
            \dot{x}(t) = {} & a x(t) + b u(t, x(t))
            \label{eqn:example-system:dynamics}
            \\
            u(t, x) = {}
            \begin{split}
               & {}\mathrel{\phantom{+}} \sqrt{\omega}\, k^{\tfrac 1 2} |x| \sin( \log(\tfrac 1 2 x^2)) \sin(\omega t) \\
               & {}+ \sqrt{\omega} \, k^{\tfrac 1 2} |x| \cos( \log(\tfrac 1 2 x^2)) \cos(\omega t)
            \end{split}
          \end{align}
        \end{subequations}
        with \( a \in \mathbb{R}\), \(b \in \mathbb{R}\!\setminus\!\set{0}\) unknown, yet with \(k\) sufficiently large, i.e.\ \( k \in \interval[open]{\tfrac{a}{b^2}}{\infty}\).
        By convention, let \(|x| \sin( \log(\frac 1 2 x^2)) |_{x=0} = |x| \cos( \log(\frac 1 2 x^2)) |_{x=0} = 0\).
  \item[Numerical Setup]
        We simulate system~\eqref{eqn:example-system} with \( {t \in \interval[]{\SI{0}{\second}}{\SI{10}{\second}}}\), \({ \omega = \SI{200}{\radian\per\second}}\), \( {a= 2}\), \({ b = -3}\) and \({ k = \tfrac 1 2}\) using an \textit{ode1}-solver with stepwidth~\({h=\SI{1e-4}{\second}}\).
  \item[Expected Results]
        System~\eqref{eqn:example-system} is in the form~\eqref{eqn:es-system} with vector fields
        \begin{subequations}
          \label{eqn:example-system:vectorfields}
          \begin{align}
            f_0(t, x) & = a x,
            \label{eqn:example-system:0}                                                                   \\
            f_1(t, x) & = bk^{\tfrac 1 2} |x| \sin( \log (\tfrac 1 2 x^2) ), 	\label{eqn:example-system:1} \\
            f_2(t, x) & = bk^{\tfrac 1 2} |x| \cos( \log (\tfrac 1 2 x^2) ),
            \label{eqn:example-system:2}
          \end{align}
          powers \( p_0 = 0\), \( p_1 = p_2 = \tfrac 1 2\) and dither signals
          \begin{multicols}{2}
            \noindent
            \begin{equation}
              u_1(s) = \sin(s),
            \end{equation}
            \begin{equation}
              u_2(s) = \cos(s)
            \end{equation}
          \end{multicols}
          \noindent
          and satisfies \Cref{ass:ui,ass:fi,ass:pi} with \( L  \leq 2838\), as shown in \Cref{lemma:lipschitz-properties-1,lemma:lipschitz-properties-2,lemma:lipschitz-properties-3}.
        \end{subequations}
        The associated Lie-bracket system is
        \begin{align}
          \bar{x}(t) & = (a-b^2 k)\bar{x}(t),
          \label{eqn:example:lbs}
        \end{align}
        which can easily be shown to be uniformly exponentially stable for our choice of values \( k \in \interval[open]{\tfrac{a}{b^2}}{\infty}\) .
        Moreover, we may select \(\bar{\alpha} = 1\) and \(\bar{\beta} = b^2 k - a\).
        Then, for any \( k \in \interval[open]{\tfrac{a}{b^2}}{\infty}\), any \( t_f \in \interval[open]{0}{\infty}\), we may choose \(D  = \tfrac 1 2 (e^{-\bar \beta t_f} -1)\)
        to satisfy~\eqref{eqn:fb1:choice-of-tf-and-d}. Accordingly, we can select a sufficiently large \( \omega^\star\) using~\eqref{eqn:def:omega-star}.

        \Cref{lemma:global-exponential-stability-result} asserts, that for all \( \omega \in \interval[open]{\omega^\star}{\infty}\), the system~\eqref{eqn:example-system} is uniformly globally exponentially stable, and we expect our simulation results to reflect this;
        We wish to emphasize that this conclusion cannot be drawn from any of the existing results in the literature.

  \item[Observations]
        The resulting trajectories are shown in \Cref{fig:simulation-result}.
        Observe, that the hull of the solution to the example system~\eqref{eqn:example-system}~(\ref{plot:exponential-convergence:ai}) is described by an exponentially decaying curve~(\ref{plot:exponential-convergence:upperbound-emprical}), which contains the curve of the associated Lie-bracket system~\eqref{eqn:example:lbs}~(\ref{plot:exponential-convergence:lbs}).
        Notably, the convergence rate of~\(x_{\omega}(t; t_0, x_0)\) seems to be proportional to \( - \bar\beta\)~(\ref{plot:exponential-convergence:upperbound-emprical}) instead of the calculated rate~\(-\beta\)~(\ref{plot:exponential-convergence:upperbound-analytical}), which indicates a faster convergence than postulated.
        The rate of convergence may be tuned further using different values of \( k\).
\end{description}

\section{Discussion and Outlook}\label{sec:discussion-and-outlook}
 The statements in our~\Cref{lemma:approximation-lemma} are very close to the statement from~\cite[Lemma~1.1]{labar2019thesis}.
The key difference in the current work lies in the fact that we do not give an absolute error bound between solutions to~\eqref{eqn:lbs} and~\eqref{eqn:es-system} in~\eqref{eqn:lemma:approximation-lemma:approximation} over finite intervals, but instead provide an error bound depending on the distance of the common initial state~\(x_0\) to the origin.
Observe, that under the restriction of \( x_0 \) to any compact set~\(\mathcal{V}\) in our result, the original statement from~\cite{labar2019thesis} may be recovered.
Observe further, that per assumptions of \Cref{lemma:approximation-lemma} the~\eqref{eqn:lbs} is bounded over finite intervals, and said bound is uniform in time \( t_0\), yet just like the error bound, linearly dependent on the distance between \(x_0\) and the origin.
Together,~\eqref{eqn:lemma:approximation-lemma:condition} and~\eqref{eqn:lemma:approximation-lemma:approximation} enable the global existence of a sufficient minimal dither frequency~\( \omega^\star \), i.e.\ the dither frequency may be chosen uniformly in \( t_0\) and globally without the need to address \(x_0\).

Scaling the error bound with the distance between the initial state~\( x_0 \) and the origin at the starting time \( t_0 \) of a finite interval allows us to provide a \textbf{global} uniform exponential stability result in~\Cref{lemma:global-exponential-stability-result}.
In particular, this is in contrast to the \textbf{semi-global} practical stability results provided in~\cite{labar2019thesis}.
To prove~\Cref{lemma:global-exponential-stability-result}, we iteratively applied \Cref{lemma:approximation-lemma} in a similar fashion to~\cite[Lemma~1.2 and Lemma~1.3]{labar2019thesis}.
Unlike~\cite{labar2019thesis}, our scaling with the initial condition \( x_0\) at the start of finite intervals from \Cref{lemma:approximation-lemma} allows for solutions of~\eqref{eqn:es-system} to converge to zero exponentially  uniformly in time for any arbitrary initial condition \(x_0\), provided that the origin of~\eqref{eqn:lbs} is globally uniformly exponentially stable and provided a sufficiently high dither frequency~\(\omega\) has been selected.

In~\Cref{lemma:global-exponential-stability-result} we have introduced a novel technique for using  Lie-bracket approximations to infer stability properties  in the sense of Lyapunov, which notably hold \textbf{globally}.
The major novelty is that previous works were concerned with practical stability only~\cite{Drr2017ExtremumSF}\cite{abdelgalil2024initializationfreeliebracketextremumseeking}, required restricting the solutions to compacts sets, resulting in semi-global properties only \cite{SUTTNER201715464}\cite{DURR20131538}\cite{Drr2017ExtremumSF}, or obtained global convergence as well, however, for a more general class of systems and through use of adaptive dither frequencies~\cite{SUTTNER2019214}.

In \Cref{sec:simulation-example} we have provided one example system, satisfying the conditions of \Cref{lemma:global-exponential-stability-result}, which demonstrates that the applicable class of systems is non-empty.
While the sufficient minimal dither frequency~\( {\omega^\star \approx \num{1e66}}\) would have been much higher according to~\eqref{eqn:def:omega-star}, the example demonstrates the highly conservative nature of our estimations throughout the trajectory approximation.
It can be seen in~\Cref{fig:simulation-result}, that in practice far lower dither frequencies suffice to drive the example system state~\(x\) exponentially to the origin.
Moreover, in practice (as may be seen in \Cref{fig:simulation-result}) it seems sufficient, to upper-bound \( \norm{x_{\omega}(t; t_0, x_0)}\) by \( \bar\alpha \, \norm{x_0} \,e^{-\bar\beta(t- t_0)}\) instead of the more conservative bound \( \bar\alpha \, \norm{x_0} \,e^{-\beta(t- t_0)}\).
It remains an open question, whether the bounds can be tightened accordingly in general.



This paper is one more step in the direction of our ultimate goal to employ Lie-bracket approximations together with high-gain adaptive control for stabilization of completely unknown systems.

However, the provided  results are useful by themselves, as they allow several vibrational control results under stricter assumptions to be extended for \textbf{global} uniform exponential stability.




 \nocite{}
 \bibliographystyle{IEEEtran}
 \bibliography{references}             

 \appendices%

 \renewcommand{\thesectiondis}[2]{\Alph{section}:}

\section{Auxiliary Results}\label{sec:auxiliary-results}
 This section contains extended auxiliary statements needed for the functional expansion.

\begin{assumption}[Assumptions for functional expansion]\label{ass:functional-expansion}
  Consider systems of form~\eqref{eqn:es-system}.
  Let \( t_0 \in \mathbb{R}\), \( t_1 \in \interval[open]{t_0}{\infty}\), \( \omega_1 \in \interval[open right]{1}{\infty}\), \( \Lambda_1 \in \interval[open right]{1}{\infty} \) and \( x_0 \in \mathbb{R}^n\), such that the trajectory \( x_{\omega_1}(t; t_0, x_0)\) of~\eqref{eqn:es-system} through \( x(t_0) = x_0\) exists on \( \interval[]{t_0}{t_1} \) and satisfies
  \begin{align}
    \norm {x_{\omega_1}(t; t_0, x_0)} & \leq \Lambda_1 \norm{x_0}\label{eqn:es-solution-bound}
  \end{align}
  for all \( t \in \interval[]{t_0}{t_1}\). Define
  \(
  T_1  \coloneq 2\pi \omega_1^{-1} \),
  \( r  \coloneq \floor*{\tfrac{t_1 - t_0}{T_1}} \),
  and
  \( t_q  \coloneq t_0+ qT_1\) ,
  for all \( q \in \set{0, 1, \ldots, r}\).    Let \( t_s \in \interval[]{t_0}{t_r}\) and \({t_e \in \interval[]{t_s}{ \min{\{t_s+T_1, t_1\}}}} \).
  Let  \( v_i: \mathbb{R} \to \mathbb{R}\), \( V_i: \mathbb{R} \times \mathbb{R} \to \mathbb{R}\) and \( V_{ij}: \mathbb{R} \times \mathbb{R} \to \mathbb{R}\) be defined by
  \begin{align}
    v_i(t)
     & \coloneq
    \omega^{p_i}\,  u_i(\omega t)
    \label{eqn:def:vi}
    \\
    V_i(t_s, t_e)
     & \coloneq
    \int_{t_s}^{t_e} v_i(\theta) \odif{\theta}
    \label{eqn:def:Vi}
    \\
    V_{ij}(t_s, t_e)
     & \coloneq
    \int_{t_s}^{t_e} v_i(\theta) \, V_j(t_s, \theta) \odif{\theta}.
    \label{eqn:def:Vij}
  \end{align}
  Moreover, define
  \begin{subequations}
    \begingroup
    \allowdisplaybreaks%
    \mathtoolsset{multlined-width=\displaywidth}
    \begin{align}
      \begin{split}
        R(t_s, t_e)
        \coloneq{} & 
        \smash[b]{
        \sum_{\mathclap{\substack{1 \leq i \leq l \\
                        0 \leq j, m  \leq l}}}
          }
        \omega_1^{p_i + p_j + p_m}
        \!\!\!
        \int_{t_s}^{t_e} \!\!\!\! u_i(\omega_1 \theta)
        \int_{t_s}^{\theta} \!\!\!\! u_j(\omega_1 \tau)
        \\
        \MoveEqLeft[4]
        \makebox[\displaywidth-4em][r]{
            \(
            \displaystyle
            \times
            \smash[b]{
                \int_{t_s}^{\tau}
              } \!\!\! \mathcal{L}_{f_m}\mathcal{L}_{f_j}f_i(\sigma, x_{\omega_1}(\sigma; t_0, x_0))
            \)
          }
        \\
        \MoveEqLeft[5]
        \makebox[\displaywidth-5em][r]{
            \(
            \displaystyle
            \times \,
            u_m(\omega_1 \sigma)
            \odif{\sigma}
            \odif{\tau}
            \odif{\theta},
            \)
          }
      \end{split}
      \label{eqn:functional-expansion:R}
      \\
      %
      %
      \begin{split}
        R_{T1}(t_0, t_1)
        \coloneq{} & 
        \sum_{q=0}^{r-1} \bigl\lbrack
        Q_{V_1}(t_q, t_{q+1})
        + Q_{V_2}(t_q, t_{q+1})
        \\[-1em]
        \MoveEqLeft[4]
        \makebox[\displaywidth-4em][r]{
            \(
            \displaystyle
            \mathrel{+} Q_1(t_q, t_{q+1})
            + R_{L1}(t_q, t_{q+1})
            \bigr\rbrack
            \)
          }
        \\
                   & 
        \mathrel{-} I(t_r, t_1)
        + Q_{v_1}(t_r, t_1)
        + Q_{V_2}(t_r, t_1)
        \\
                   & 
        + Q_1(t_r, t_1)
        + Q_{1T}(t_r, t_1)
        + H(t_r, t_1),
      \end{split}
      \label{eqn:functional-expansion:RT1}
      \\
      \begin{split}
        Q_{V_1}(t_s, t_e)
        \coloneq{} & 
        \smash[b]{
            \sum_{i=1}^l \omega_1^{p_i}\! \int_{t_s}^{t_e} \! \int_{t_s}^{\theta}
          }
        \pdif{1} f_i(\tau, x_{\omega_1}(\tau; t_0, x_0))
        \\
        \MoveEqLeft[4]
        \makebox[\displaywidth-4em][r]{
            \(
            \displaystyle
            \mathrel{\times} u_i(\omega_1 \theta) \odif{\tau} \odif{\theta}
            \)
          }
      \end{split}
      \label{eqn:functional-expansion:QV1}
      \\
      \begin{split}
        Q_{V_2}(t_s, t_e)
        \coloneq{} & 
        \sum_{\mathclap{\begin{limitarray}
                                  1&\leq i &< l \\
                                  0&\leq j &\leq l
                                \end{limitarray}}}
        \omega_1^{p_i + p_j}
        \int_{t_s}^{t_e} u_i(\omega_1 \theta)
        \int_{t_s}^\theta u_j(\omega_1 \tau)
        \\
        \MoveEqLeft[4]
        \makebox[\displaywidth-4em][r]{
            \(
            \displaystyle
            \mathrel{\times}
            \smash[b]{
                \int_{t_s}^\tau
              } \pdif{1} \mathcal{L}_{f_j} f_i(\sigma, x_{\omega_1}(\sigma; t_0, x_0))
            \odif{\sigma}
            \odif{\tau}
            \odif{\theta}
            \)
          }
      \end{split}
      \label{eqn:functional-expansion:QV2}
      \\
      %
      %
      %
      \begin{split}
        Q_1(t_s, t_e)
        \coloneq{} & 
        \smash{\sum_{i=1}^l }\mathcal{L}_{f_0}f_i(t_s, x_{\omega_1}(t_s; t_0, x_0))
        \mathrel{\times} V_{i0}(t_s, t_e)
      \end{split}
      \label{eqn:functional-expansion:Q1}
      \\
      \begin{split}
        Q_{1T}(t_s, t_e)
        \coloneq{} & 
        \smash[b]{
            \sum_{\mathclap{\begin{limitarray}
                                      1& \leq i &< l \\
                                      i &< j &\leq l
                                    \end{limitarray} }}
          }
        \mathcal{L}_{f_j}f_i(t_s, x_{\omega_1}\!(t_s; t_0, x_0))
        \\
        \MoveEqLeft[4]
        \makebox[\displaywidth-4em][r]{
            \(
            \displaystyle
            \mathrel{\times} V_i(t_s, t_e)\, V_j(t_s, t_e)
            \)
          }
        \\
                   & \mathrel{\phantom{\coloneq}}\! \mathllap{+\: \tfrac 1 2}
        \smash[b]{\sum_{i=1}^l} \mathcal{L}_{f_i} f_i(t_s, x_{\omega_1}\!(t_s; t_0, x_0))
        \,
        V_i^2(t_s, t_e)
      \end{split}
      \label{eqn:functional-expansion:Q1T}
      \\
      \begin{split}%
        H(t_s, t_e)
        \coloneq{} & 
        \sum_{i=1}^l f_i(t_s, x_{\omega_1}(t_s; t_0, x_0))\, V_i(t_s, t_e)
        \\*
                   & +  \smash[b]{
            \sum_{\mathclap{\begin{limitarray}
                                      1 &\leq i &< l \\
                                      i &< j &\leq l
                                    \end{limitarray}}}
          }
        \commutator{f_i}{f_j}(t_s, x_{\omega_1}(t_s; t_0, x_0))
        \,
        V_{ji}(t_s, t_e)
      \end{split}%
      \label{eqn:functional-expansion:H}
      \\
      \begin{split}
        I(t_s, t_e)
        \coloneq{} & 
        \smash[b]{
            \!
            \int_{t_s}^{t_e}
            \!
            \sum_{\mathclap{\begin{limitarray}
                                      1 & \leq i &< l\\
                                      i &< j &\leq l
                                    \end{limitarray}}}
          }
        \commutator{f_i}{ f_j}(\theta, x_{\omega_1}(\theta; t_0, x_0))
        \,
        \gamma_{ij}(\omega_1) \odif{\theta}
      \end{split}
      \label{eqn:functional-expansion:I}
      \\
      \intertext{and}
      \begin{split}
        R_{L1}(t_s, t_e)
        \coloneq{} & 
        \smash[b]{
            \sum_{\mathclap{\begin{limitarray}
                                      1 &\leq i &< l \\ i &< j &\leq l
                                    \end{limitarray}}}
          }
        \gamma_{ij}(\omega_1)
        \!\!
        \int_{t_s}^{t_e}
        \Bigl\lparen
        \commutator{f_i}{f_j}(t_s, x_{\omega_1}(t_s; t_0, x_0))
        \\
        \MoveEqLeft[4]
        \makebox[\displaywidth-4em][r]{
        \(
        \displaystyle
        \mathrel{-}
        \commutator{f_i}{f_j}(\theta, x_{\omega_1}(\theta; t_0, x_0))
          \smash[b]{\Bigr\rparen}
          \odif{\theta}
          \)
          }.
      \end{split}
      \label{eqn:functional-expansion:RL1}
    \end{align}
    \endgroup
  \end{subequations}
\end{assumption}

\begin{lemma}[Functional Expansion of~\eqref{eqn:es-system}]\label{lemma:functional-expansion}~
  Consider systems of form~\eqref{eqn:es-system} and
  suppose, that \Cref{ass:functional-expansion} holds.
  Then, the trajectory \( x_{\omega_1}(t; t_0, x_0)\) may be expressed as
  \begin{align}
    \begin{split}
      \MoveEqLeft
      x_{\omega_1}(t_1; t_0, x_0) = x_0  +
      \int_{t_0}^{t_1}
      f_0(\theta, x_{\omega_1}(\theta; t_0, x_0))
      \\
       & 
      +
      \sum_{\mathclap{\begin{limitarray}
                                1 & \leq i & < l \\
                                i & < j     & \leq l
                              \end{limitarray}}}
      \commutator{f_i}{f_j}(\theta, x_{\omega_1}(\theta; t_0, x_0))\,\gamma_{ij}(\omega_1) \odif{\theta}
      \\[-0.25em]
       & 
      +
      \sum_{q=0}^{r-1} R(t_q, t_{q+1})
    \end{split}
    + R(t_r, t_1) + R_{T_1}(t_0, t_1)
    \label{eqn:functional-expansion:main}
  \end{align}
  with \(\gamma_{ij}(\omega_1)\) defined in~\eqref{eqn:gamma-ij} for \( \omega=\omega_1\).
\end{lemma}
The proof has been omitted in favour of the page limit. The complete derivation is given in \cite{labar2019thesis}.

\begin{theoremEnd}[restate]{lemma}[Auxiliary bounds]\label{prop:auxiliary-bounds}
  Let \( i, j \in \{0, 1, \ldots, l\}\) and let \(u_i\) satisfy \Cref{ass:ui}. Consider the functions \( v_i: \mathbb{R} \to \mathbb{R}\), \( V_i: \mathbb{R} \times \mathbb{R} \to \mathbb{R}\) and \( V_{ij}: \mathbb{R} \times \mathbb{R} \to \mathbb{R}\) as defined in~\Cref{ass:functional-expansion}
  \begin{subequations}
    Then,
    \begin{align}
      \abs{ v_i(t)}
      & \leq
      \phantom{\tfrac{1 }{2}}
      \omega^{p_i}
      \label{eqn:bound:vi}
      \\
      \abs{V_i(t_s, t_e)}
      & \leq
      \phantom{\tfrac{1 }{2}}
      \omega_1^{p_i} {(t_e -t_s)}
      \label{eqn:bound:Vi}
      \\
      \abs{ V_{ij}(t_s, t_e)}
      & \leq
      \tfrac{1 }{2}
      \omega^{p_i + p_j }{(t_e - t_s)}^2.
      \label{eqn:bound:Vij}
    \end{align}
  \end{subequations}
\end{theoremEnd}
\begin{theoremEnd}[%
    category={auxiliary-result},%
  ]{lemma}[Bounds on \(Q_{V_1}\)]\label{lemma:bounds:QV1}
  Consider \(Q_{V_1}\)  as defined in~\eqref{eqn:functional-expansion:QV1} with \Cref{ass:functional-expansion}. Then,
  \begin{align}
    \norm{Q_{V_1}(t_s, t_e)}
    & \leq
    \pi l  L \Lambda_1 \, \norm{  x_0  } \, \omega_1^{p'-1} \, T_1.
    \qedhere
  \end{align}
\end{theoremEnd}
\begin{proofEnd}
  The result is established by straightforwardly calculating
  \begin{align*}
    \norm{ Q_{V_1}(t_s, t_e) }
    & \labelrel{=}{rel:tmp:QV1:1}
    \sum_{i=1}^l
    \begin{multlined}[t]
      \omega_1^{p_i} \int_{t_s}^{t_e} \int_{t_s}^{\theta} \norm{ \pdif{1} f_i(\tau, x_{\omega_1}(\tau; t_0, x_0)) }
      \\
      \times \abs{ u_i(\omega_1 \theta) }  \odif{\tau} \odif{\theta}
    \end{multlined}
    \\
    & \labelrel{=}{rel:tmp:QV1:2}
    l {\omega_1}^{p'} \int_{t_s}^{t_e} \int_{t_s}^{\theta}  L \, \norm{  x_{\omega_1}(\tau; t_0, x_0) }
    \odif{\tau} \odif{\theta}
    \\
    & \labelrel{=}{rel:tmp:QV1:3}
    \tfrac 1 2 l   L \Lambda_1 \, \norm{  x_0  } \, \omega_1^{p'}{(t_e - t_s)}^2
    \\
    & \labelrel{=}{rel:tmp:QV1:4}
    \tfrac 1 2 l   L \Lambda_1 \, \norm{  x_0  } \, \omega_1^{p'} T_1^2
    \\
    & \labelrel{=}{rel:tmp:QV1:5}
    \pi l  L \Lambda_1  \, \norm{  x_0  } \, \omega_1^{p'-1} \, T_1
  \end{align*}
  where we
  \begin{steps}
    \item[\eqref{rel:tmp:QV1:1}] use absolute homogenicity and triangle inequality of norms,
    \item[\eqref{rel:tmp:QV1:2}] substitute in~\eqref{eqn:def:p:1}, use \Cref{ass:fi:2:4} and \Cref{ass:ui:1},
    \item[\eqref{rel:tmp:QV1:3}] apply~\eqref{eqn:es-solution-bound} and solve the integral,
    \item[\eqref{rel:tmp:QV1:4}] remember that \( t_e - t_s \leq T_1\) by the assumptions of this lemma,
    \item[\eqref{rel:tmp:QV1:5}] and finally remember the definition of \( T_1\),
  \end{steps}
  which concludes this lemma.
\end{proofEnd}
\begin{theoremEnd}[%
    category=auxiliary-result,%
  ]{lemma}[Bounds on \(Q_{V_2}\)]\label{lemma:bounds:QV2}
  Consider \(Q_{V_2}\)  as defined in~\eqref{eqn:functional-expansion:QV2} with~\Cref{ass:functional-expansion}. Then,
  \begin{align}
    \norm{Q_{V_2}(t_s, t_e)}
    & \leq
    \tfrac{2}{3}
    \pi^2
    l^2
    L
    \Lambda_1 \, \norm{x_0}
    \,
    \omega_1^{p'-1}
    \,
    {T_1}.
  \end{align}
\end{theoremEnd}
\begin{proofEnd}
  The result is obtained by straightforwardly calculating
  \begin{align*}
    \norm{ Q_{V_2}(t_s, t_e) }
    & \labelrel{\leq}{rel:tmp:QV2:1}
    \sum_{\mathclap{\begin{limitarray}
              1&\leq i &< l \\
              0&\leq j &\leq l
            \end{limitarray}}}
    \begin{multlined}[t]
      \omega_1^{p_i + p_j}
      \int_{t_s}^{t_e} \abs{u_i(\omega_1 \theta) }
      \int_{t_s}^\theta \abs{ u_j(\omega_1 \tau) }
      \\
      \times \int_{t_s}^\tau \norm{ \pdif{1}\mathcal{L}_{f_j} f_i(\sigma, x_{\omega_1}(\sigma; t_0, x_0)) }
      \\
      \dots
      \odif{\sigma}
      \odif{\tau}
      \odif{\theta}
    \end{multlined}
    \\
    & \labelrel{\leq}{rel:tmp:QV2:2}
    \begin{multlined}[t]
      (l-1)(l+1)
      L
      \omega_1^{p''}
      \\
      \times
      \int_{t_s}^{t_e} \!
      \int_{t_s}^\theta
      \int_{t_s}^\tau \!\! \norm{ x_{\omega_1}(\sigma; t_0, x_0) }
      \odif{\sigma}
      \odif{\tau}
      \odif{\theta}
    \end{multlined}
    \\
    & \labelrel{\leq}{rel:tmp:QV2:3}
    \tfrac{1}{6}
    l^2
    L
    \Lambda_1\, \norm{x_0} \,
    \omega_1^{p''}
      {(t_e - t_s)}^3
    \\
    & \labelrel{\leq}{rel:tmp:QV2:4}
    \tfrac{1}{6}
    l^2
    L
    \Lambda_1 \, \norm{x_0} \,
    \omega_1^{p''}
      {T_1}^3
    \\
    & \labelrel{\leq}{rel:tmp:QV2:5}
    \tfrac{2}{3}
    \pi^2
    l^2
    L
    \Lambda_1 \, \norm{x_0} \,
    \omega_1^{p''-2}
    \,
    {T_1}
    \\
    & \labelrel{\leq}{rel:tmp:QV2:6}
    \tfrac{2}{3}
    \pi^2
    l^2
    L
    \Lambda_1 \, \norm{x_0} \,
    \omega_1^{p'-1}
    \,
    {T_1},
  \end{align*}
  where we
  \begin{steps}
    \item[\eqref{rel:tmp:QV2:1}] use absolute homogenicity and triangle inequality of norms,
    \item[\eqref{rel:tmp:QV2:2}] substitute in~\eqref{eqn:def:p:1}, used \Cref{ass:fi:2:5} and \Cref{ass:ui:1},
    \item[\eqref{rel:tmp:QV2:3}] apply~\eqref{eqn:es-solution-bound} and solve the integral, as well as overestimate the contribution of \((l-1)(l+1)\),
    \item[\eqref{rel:tmp:QV2:4}] remember that \( t_e -t_s \leq T_1\) from the assumptions,
    \item[\eqref{rel:tmp:QV2:5}] remember the definition of \( T_1\) and
    \item[\eqref{rel:tmp:QV2:6}] use~\eqref{eqn:def:p:1}, which concludes the statement of this lemma.\qedhere
  \end{steps}
\end{proofEnd}
\begin{theoremEnd}[%
    category=auxiliary-result,%
  ]{lemma}[Bounds on \(Q_{1T}\)]\label{lemma:bounds:Q1T}
  Consider \(Q_{1T}\)  as defined in~\eqref{eqn:functional-expansion:Q1T} with \Cref{ass:functional-expansion}. Then,
  \begin{align}
    \norm{ Q_{1T}(t_s, t_e) }
    & \leq
    2 \pi^2  l^2
    L \Lambda_1 \, \norm{x_0} \,
    \omega_1^{p'-1}.
  \end{align}
\end{theoremEnd}
\begin{proofEnd}
  The result is established by straightforwardly calculating
  \begingroup
  \allowdisplaybreaks
  \begin{align*}
    \norm{ Q_{1T}(t_s, t_e) }
    %
    \labelrel{\leq}{rel:tmp:Q1T:1}
      {}
    \begin{split}
      &
      \smash[b]{
          \sum_{\mathclap{\begin{limitarray}
                    1 &\leq i &< l \\
                    i &< j &\leq l
                  \end{limitarray}}}
        }
      \norm{ \mathcal{L}_{f_j}f_i(t_s, x_{\omega_1}(t_s; t_0, x_0)) } \,
      \\
      \MoveEqLeft[4]
      \makebox[\displaywidth-4em][r]{
          \(
          \displaystyle
          \mathrel{\times} | V_i(t_s, t_e)| \,| V_j(t_s, t_e) |
          \)
        }
      \\
      & {+}
      \tfrac 1 2
      \smash[b]{
          \sum_{\mathclap{i=1}}^l
        }
      \norm{ \mathcal{L}_{f_i\vphantom{f_j}} f_i(t_s, x_{\omega_1}(t_s; t_0, x_0)) }
      \\
      \MoveEqLeft[4]
      \makebox[\displaywidth-4em][r]{
          \(
          \displaystyle
          \mathrel{\times}
            {| V_i(t_s, t_e) |}^2
          \)
        }
    \end{split}
    \\
    %
    \labelrel{\leq}{rel:tmp:Q1T:2}
      {}
    \begin{split}
      &
      \tfrac 1 2 l(l-1)
      L\, \norm{x_{\omega_1}(t_s; t_0, x_0)}
      \,
      \omega_1^{2p'} {(t_e -t_s)}^2
      \\
      &+
      \tfrac 1 2 l
      L \, \norm{x_{\omega_1}(t_s; t_0, x_0)}
      \,\omega_1^{2p'} {(t_e -t_s)}^2
    \end{split}
    \\
    %
    \labelrel{\leq}{rel:tmp:Q1T:3}
      {} &
    \tfrac 1 2 l^2
    L \, \norm{x_{\omega_1}(t_s; t_0, x_0)}
    \, \omega_1^{2p'}\!\! {(t_e -t_s)}^2
    \\
    %
    \labelrel{\leq}{rel:tmp:Q1T:4}
      {} &
    \tfrac 1 2 l^2
    L \, \norm{x_{\omega_1}(t_s; t_0, x_0)}
    \, \omega_1^{2p'} T_1^2
    \\
    %
    \labelrel{\leq}{rel:tmp:Q1T:5}
      {} &
    2 \pi^2  l^2
    L \, \norm{x_{\omega_1}(t_s; t_0, x_0)}
    \, \omega_1^{2p'-2}
    \\
    \labelrel{\leq}{rel:tmp:Q1T:6}
      {} &
    2 \pi^2  l^2
    L \Lambda_1 \, \norm{x_0}
    \, \omega_1^{p'-1},
  \end{align*}
  \endgroup
  where we
  \begin{steps}
    \item[\eqref{rel:tmp:Q1T:1}] use submultiplicativity and triangle equation,
    \item[\eqref{rel:tmp:Q1T:2}] insert \Cref{ass:fi:2:2}, \eqref{eqn:def:p:1} and~\eqref{eqn:bound:Vi} from \Cref{prop:auxiliary-bounds},
    \item[\eqref{rel:tmp:Q1T:3}]  sum the terms and conservatively overestimate the contribution of \(l\),
    \item[\eqref{rel:tmp:Q1T:4}]  remember \( t_e - t_s \leq T_1\) from the assumptions of this lemma,
    \item[\eqref{rel:tmp:Q1T:5}] the definition of \(T_1\) and finally
    \item[\eqref{rel:tmp:Q1T:6}]  use~\eqref{eqn:es-solution-bound} and \( p' -1 \leq 0\), which concludes this lemma.\qedhere
  \end{steps}
\end{proofEnd}
\begin{theoremEnd}[%
    category=auxiliary-result,%
  ]{lemma}[Bounds on \( Q_1 \)]\label{lemma:bounds:Q1}
  Consider \(Q_{1}\)  as defined in~\eqref{eqn:functional-expansion:Q1} with \Cref{ass:functional-expansion}. Then,
  \begin{align}
    \norm{Q_1(t_s, t_e)}
    & \leq
    \pi l {L} \Lambda_1 \, \norm{x(t_0)} \,
    \omega_1^{p'-1} \, T_1.
  \end{align}
\end{theoremEnd}
\begin{proofEnd}
  The result is established by straighforwardly calculating
  \begingroup
  \allowdisplaybreaks
  \begin{align*}
    \norm{Q_1(t_s, t_e)}
    %
    & \labelrel{\leq}{rel:tmp:Q1:1}
    \sum_{i=1}^l \norm{ \mathcal{L}_{f_0}f_i(t_s, x_{\omega_1}(t_s; t_0, x_0))}  \, | V_{i0}(t_s, t_e)|
    \\
    %
    & \labelrel{\leq}{rel:tmp:Q1:2}
    \tfrac 1 2
    \smash[b]{
        \sum_{i=1}^l
      }
    \norm{ \mathcal{L}_{f_0}f_i(t_s, x_{\omega_1}(t_s; t_0, x_0))}
    \\
    \MoveEqLeft[3]
    \makebox[\displaywidth-3em][r]{
        \(
        \displaystyle
        \smash[t]{
            \mathrel{\times}
            \omega_1^{p'}
              {(t_e -t_s)}^2
          }
        \)
      }
    \\
    %
    & \labelrel{\leq}{rel:tmp:Q1:3}
    \tfrac 1 2 l {L} \,  \norm{x_{\omega_1}(t_s; t_0, x_0)}
    \, \omega_1^{p'} {(t_e -t_s)}^2
    \\
    %
    & \labelrel{\leq}{rel:tmp:Q1:4}
    \pi l {L} \, \norm{x_{\omega_1}(t_s; t_0, x_0)}
    \, \omega_1^{p'-1} \, T_1
    \\
    & \labelrel{\leq}{rel:tmp:Q1:5}
    \pi l {L} \Lambda_1 \, \norm{x_0} \,
    \omega_1^{p'-1} \, T_1,
  \end{align*}
  \endgroup
  where we use
  \begin{steps}
    \item[\eqref{rel:tmp:Q1:1}] triangle equality and submultiplicativity of norms,
    \item[\eqref{rel:tmp:Q1:2}] property~\eqref{eqn:bound:Vij} from \Cref{prop:auxiliary-bounds} in conjunction with \(p_0=0\) and \Cref{eqn:def:p:1} ,
    \item[\eqref{rel:tmp:Q1:3}] \Cref{ass:fi:2:3},
    \item[\eqref{rel:tmp:Q1:4}] remember \(t_e -t_s \leq T_1\) from the assumptions of this lemma and the definition of \(T_1\),
    \item[\eqref{rel:tmp:Q1:5}] use~\eqref{eqn:es-solution-bound}, which concludes this lemma.\qedhere
  \end{steps}
\end{proofEnd}

\begin{theoremEnd}[%
    category=auxiliary-result,%
  ]{lemma}[Bounds on \(H\)]\label{lemma:bounds:H} {}
  Consider \(H\)  as defined in~\eqref{eqn:functional-expansion:H} with \Cref{ass:functional-expansion}. Then,
  \begin{align}
    \norm{ H(t_s, t_e) }
    & \leq
    8 \pi^2 l^2 L \Lambda_1 \, \norm{x_0} \,  \omega_1^{p'-1}.
  \end{align}
\end{theoremEnd}
\begin{proofEnd}
  The result is established by straightforwardly calculating
  \begingroup
  \allowdisplaybreaks
  \begin{align*}
    \norm{ H(t_s, t_e) }
    %
    \labelrel{\leq}{rel:tmp:H:1}
      {}
    \begin{split}
      &
      \sum_{i=1}^l
      \norm{ f_i(t_s, x_{\omega_1}(t_s; t_0, x_0)) }
      \,
      {| V_i(t_s, t_e) |}
      \\
      & +
      \smash[b]{
          \sum_{\mathclap{\begin{limitarray}
                    1 &\leq i &< l \\
                    i &< j &\leq l
                  \end{limitarray}}}
        }
      \norm{ \commutator{f_i}{f_j}(t_s, x_{\omega_1}(t_s; t_0, x_0)) }
      \\
      \MoveEqLeft[3]
      \makebox[\displaywidth-3em][r]{
          \(
          \displaystyle
          \mathrel{\times}
            {| V_{ji}(t_s, t_e) |}
          \)
        }
    \end{split}
    \\
    %
    \labelrel{\leq}{rel:tmp:H:2}
      {}
    \begin{split}
      &
      \sum_{i=1}^l
      \begin{multlined}[t]
        L \,\norm{x_{\omega_1}(t_s; t_0, x_0)}
        \,  \omega_1^{p_i} {(t_e -t_s)}
      \end{multlined}
      \\
      & +
      2 \sum_{\mathclap{\begin{limitarray}
                1 &\leq i &< l \\
                i &< j &\leq l
              \end{limitarray}}}
      \begin{multlined}[t]
        L \, \norm{ x_{\omega_1}(t_s; t_0, x_0) }
        \,   \omega_1^{p_i + p_j} {(t_e -t_s)}^2
      \end{multlined}
    \end{split}
    \\
    %
    \labelrel{\leq}{rel:tmp:H:3}
      {}
    \begin{split}
      &
      l^2 L \, \norm{x_{\omega_1}(t_s; t_0, x_0)}  \, \omega_1^{p'} {(t_e -t_s)}
      \\
      &
      +
      l^2 L \, \norm{ x_{\omega_1}(t_s; t_0, x_0) } \,  \omega_1^{2p'} {(t_e -t_s)}^2
    \end{split}
    \\
    %
    \labelrel{\leq}{rel:tmp:H:4}
      {}
    \begin{split}
      &
      2 \pi l^2 L \, \norm{x_{\omega_1}(t_s; t_0, x_0)}  \,  \omega_1^{p'-1}
      \\
      & +
      {(2\pi)}^2 l^2 L \,  \norm{ x_{\omega_1}(t_s; t_0, x_0) }   \, \omega_1^{p'-1}
    \end{split}
    \\
    %
    \leq
      {}
    \begin{split}
      &
      {(2 \pi)}^2 l^2 L \,
      \norm{x_{\omega_1}(t_s; t_0, x_0)}
      \,    \omega_1^{p'-1}
      \\
      & +
      {(2 \pi)}^2 l^2 L \,
      \norm{ x_{\omega_1}(t_s; t_0, x_0) }
      \,   \omega_1^{p'-1}
    \end{split}
    \\
    %
    \labelrel{\leq}{rel:tmp:H:5}
      {}
    &
    8 \pi^2 l^2 L \Lambda_1  \, \norm{x_0} \,   \omega_1^{p'-1},
  \end{align*}
  \endgroup
  where we
  \begin{steps}
    \item[\eqref{rel:tmp:H:1}] use triangle inequality of norms and submultiplicativity,
    \item[\eqref{rel:tmp:H:2}] insert~\eqref{eqn:bound:Vi} and~\eqref{eqn:bound:Vij} from \Cref{prop:auxiliary-bounds} as well as use the Lipschitz properties of \Cref{ass:fi:2:1,ass:fi:2:2},
    \item[\eqref{rel:tmp:H:3}] overestimate using~\eqref{eqn:def:p:1}, \( {p_i + p_j < 2 p'}\) for all \( i, i \in \set{0, 1, \ldots, l}^2\) and \(l \leq l^2\), \( l(l-1) \leq l^2\),
    \item[\eqref{rel:tmp:H:4}]  remember \( t_e - t_s \leq T_1\) from the assumptions of this lemma and use the definition of \(T_1\), remember \(\omega_1 \in \interval[open right]{1}{\infty}\) and \( p' \in \interval[]{0}{1}\),
    \item[\eqref{rel:tmp:H:5}] sum up the terms and use~\eqref{eqn:es-solution-bound}, which concludes the statement of this lemma. \qedhere
  \end{steps}
\end{proofEnd}
\begin{theoremEnd}[%
    category=auxiliary-result,%
  ]{lemma}[Bounds on \(I\)]\label{lemma:bounds:I}
  Consider \(I\)  as defined in~\eqref{eqn:functional-expansion:I} with \Cref{ass:functional-expansion}.
  Then,
  \begin{align}
    \norm{ I(t_s, t_e) }
    & \leq
    2 \pi^2  l^2  L \Lambda_1 \, \norm{  x_0 } \,   \omega_1^{p'-1}.
  \end{align}
\end{theoremEnd}
\begin{proofEnd}
  The result is established by straighforwardly calculating
  \begingroup
  \allowdisplaybreaks
  \begin{align*}
    \norm{I(t_s, t_e)}
    & =
    \norm{
    \int_{t_s}^{t_e} \sum_{\mathclap{\begin{limitarray}
              1 &\leq i &< l \\
              i &< j &\leq l
            \end{limitarray}}} \commutator{f_i}{f_j}(\theta, x_{\omega_1}(\theta; t_0, x_0)) \,\gamma_{ij}(\omega_1) \odif{\theta}
      }
    \\
    %
    & \labelrel{\leq}{rel:tmp:I:1}
    \int_{t_s}^{t_e}
    \!
    \sum_{\mathclap{\begin{limitarray}
              1 &\leq i &< l \\
              i &< j &\leq l
            \end{limitarray}}}
    \!
    \norm{ \commutator{f_i}{f_j}(\theta, x_{\omega_1}(\theta; t_0, x_0)) }
    \,
    \abs{\gamma_{ij}(\omega_1) }  \odif{\theta}
    \\
    %
    & \labelrel{\leq}{rel:tmp:I:2}
    \int_{t_s}^{t_e}
    \!  \sum_{\mathclap{\begin{limitarray}
              1 &\leq i &< l \\
              i &< j &\leq l
            \end{limitarray}}}
    L \,
    \norm{ x_{\omega_1}(\theta; t_0, x_0) }  \,
    \omega_1^{p_i + p_j} {T_1}  \odif{\theta}
    \\
    %
    & \labelrel{\leq}{rel:tmp:I:3}
    \sum_{\mathclap{\begin{limitarray}
              1 &\leq i &< l \\
              i &< j &\leq l
            \end{limitarray}}}
    L
    \Lambda_1
    \,
    \norm{ x_0 }
    \omega_1^{2p'} {T_1}
    \int_{t_s}^{t_e}  \odif{\theta}
    \\
    %
    & \labelrel{\leq}{rel:tmp:I:4}
    2 \pi^2  l^2  L \Lambda_1 \, \norm{  x_0 }     \omega_1^{2p'-2}
    \\
    %
    & \labelrel{\leq}{rel:tmp:I:5}
    2 \pi^2  l^2  L \Lambda_1 \, \norm{  x_0 }  \,   \omega_1^{p'-1},
  \end{align*}
  \endgroup
  where we
  \begin{steps}
    \item[\eqref{rel:tmp:I:1}] use absolute homogenicity and triangle equation of norms,
    \item[\eqref{rel:tmp:I:2}] insert~\eqref{eqn:gamma-ij}, and use~\cref{ass:fi:2:2} from \Cref{ass:fi},
    \item[\eqref{rel:tmp:I:3}] use~\eqref{eqn:es-solution-bound} and~\eqref{eqn:def:p:1},
    \item[\eqref{rel:tmp:I:4}] solve the integral, remember \( t_e - t_s \leq T_1 \) from the assumptions of this lemma, use the definition of \( T_1\) and express everything in terms of \( \omega_1\),
    \item[\eqref{rel:tmp:I:5}] and finally remember \( p' \in \interval[open]{0}{1}\) from \Cref{ass:pi}, which concludes the statement of this lemma. \qedhere
  \end{steps}
\end{proofEnd}
\begin{theoremEnd}[%
    category=auxiliary-result,%
  ]{lemma}[Bounds on \(R_{L1}\)]\label{lemma:bounds:RL1}
  Consider \(R_{L1}\)  as defined in~\eqref{eqn:functional-expansion:RL1} with \Cref{ass:functional-expansion}.
  Then,
  \begin{align}
    \norm{ R_{L1}(t_s, t_e) }
    & \leq
    \pi^2
    {(l+1)}^2
    L^2 \Lambda_1 \, \norm{ x_0 } \,
    \omega_1^{p'-1} \,
    T_1.
  \end{align}
\end{theoremEnd}
\begin{proofEnd}
  The result is established by straightforwardly calculating
  \begingroup
  \allowdisplaybreaks
  \begin{align*}
    \norm{ R_{L1}(t_s, t_e) }
    %
    & \labelrel{\leq}{rel:tmp:RL1:1}
    \smash[b]{
        \sum_{\mathclap{\begin{limitarray}
                  1 &\leq i &< l \\
                  i &< j &\leq l
                \end{limitarray}}}
      }
    \abs{ \gamma_{ij}(\omega_1) } \\*
    \MoveEqLeft[3]
    \makebox[\displaywidth-3em][r]{
    \(
    \displaystyle
    \mathrel{\times}
    \int_{t_s}^{t_e} \bigl\|
    \commutator{f_i}{f_j}(t_s, x_{\omega_1}(t_s; t_0, x_0))
      \)
      }
    \\*
    \MoveEqLeft[3]
    \makebox[\displaywidth-3em][r]{
    \(
    \displaystyle
    \mathrel{-}
    \commutator{f_i}{f_j}(\theta, x_{\omega_1}(\theta; t_0, x_0))
      \bigr\| \odif{\theta}
      \)
      }
    \\
    %
    %
    & \labelrel{\leq}{rel:tmp:RL1:3}
    \smash[b]{
        \sum_{\mathclap{\begin{limitarray}
                  1 &\leq i &< l \\
                  i &< j &\leq l
                \end{limitarray}}}
        \pi \omega_1^{0}
        \int_{t_s}^{t_e}
      }
    2L \, \bigl\| x_{\omega_1}(t_s; t_0, x_0) \\
    \MoveEqLeft[3]
    \makebox[\displaywidth-3em][r]{
        \(
        \displaystyle
        \mathrel{-}
        x_{\omega_1}(\theta; t_0, x_0) \bigr\|
        \odif{\theta}
        \)
      }
    \\
    %
    & \labelrel{\leq}{rel:tmp:RL1:4}
    2\pi
    l^2 (l+1)
    L^2 \Lambda_1 \, \norm{ x_0 }
    \omega_1^{p'}
    \int_{t_s}^{t_e} (\theta -t_s) \odif{\theta}
    \\
    %
    %
    & \labelrel{\leq}{rel:tmp:RL1:5}
    \pi
      {(l+1)}^3
    L^2 \Lambda_1 \,  \norm{ x_0 }
    \,
    \omega_1^{p'}
      {(t_e - t_s)}^2
    \\
    %
    & \labelrel{\leq}{rel:tmp:RL1:6}
    2\pi^2
    {(l+1)}^2
    L^2 \Lambda_1
    \,
    \norm{ x_0 }
    \,
    \omega_1^{p'-1}
    \,
    T_1
  \end{align*}
  \endgroup
  where we
  \begin{steps}
    \item[\eqref{rel:tmp:RL1:1}] use absolute homogenicity and the triangle inequality,
    \item[\eqref{rel:tmp:RL1:3}] remember \Cref{ass:pi}, i.e.\ that for \( p_i + p_j > 1\), either \(\gamma_{ij} = 0\) or \( \commutator{f_i}{f_j}(t, x) = 0\) for all \( t \in \mathbb{R}\) and all \( x \in \mathbb{R}^n\), insert \Cref{ass:fi:2:2} and substitute bounds of \(\gamma_{ij}\),
    \item[\eqref{rel:tmp:RL1:4}] use~\eqref{eqn:es-solution-bound} and the fundamental theorem of calculus to expand solutions of~\eqref{eqn:es-system} once,
    \item[\eqref{rel:tmp:RL1:5}] solve the integral and overestimated the contribution of \(l\),
    \item[\eqref{rel:tmp:RL1:6}] remember \( t_e -t_s \leq T_1\) from the assumptions of this lemma and use definition of \(T_1\)        to conclude the statement of this lemma. \qedhere
  \end{steps}
\end{proofEnd}
\begin{theoremEnd}[%
    category=auxiliary-result,%
  ]{lemma}[Bounds on \(R_{T1}\)]\label{lemma:bounds:RT1}
  Consider \(R_{T1}\)  as defined in~\eqref{eqn:functional-expansion:RT1} with \Cref{ass:functional-expansion}
  Then,
  \begin{align}
    \norm{ R_{T1}(t_0, t_1) }
    \leq {}
    \begin{split}
      & \pi^2 {(l+1)}^2 L^2 \Lambda_1 \,  \omega_1^{p'-1} \,  \norm{ x_0 }
      \\
      & \hphantom{\leq} \times
      \left\lparen
      8(t_1 -t_0)
      +12
      +6\pi
      \right\rparen.
    \end{split}
  \end{align}
\end{theoremEnd}
\begin{proofEnd}
  The result is established by straightforwardly calculating
  \begingroup
  \allowdisplaybreaks
  \begin{align*}
    \norm{ R_{T1}(t_s, t_e) }
    \labelrel{\leq}{rel:tmp:RT1:1}
      {}
    & \sum_{q=0}^{r-1} \bigl(
    \begin{aligned}[t]
      & \,\;\phantom{+} \norm{ Q_{V_1}(t_q, t_{q+1}) }
      + \norm{ Q_{V_2}(t_q, t_{q+1}) }
      \\
      & + \norm{Q_1(t_q, t_{q+1}) }
      + \norm{ R_{L1}(t_q, t_{q+1}) }
      \\
      & + \norm{ I(t_r, t_1) }
      + \norm{Q_{V_1}(t_r, t_1) }
      \\
      & + \norm{ Q_{V_2}(t_r, t_1) }
      + \norm{ Q_1(t_r, t_1) }
      \\
      & + \norm{ Q_{1T}(t_r, t_1) }
      +  \norm{ H(t_r, t_1) } \quad \bigr)
    \end{aligned}
    \\
    \labelrel{\leq}{rel:tmp:RT1:2}
      {}
    & \sum_{q=0}^{r-1} \bigl\lparen
    \begin{aligned}[t]
      & \pi l L \Lambda_1 \, \norm{x_0} \,   \omega_1^{p'-1} \,T_1                   \\
      & + \tfrac 2 3 \pi^2 l^2 L \Lambda_1 \, \norm{x_0} \,   \omega_1^{p'-1} \, T_1 \\
      & + \pi l L \Lambda_1 \, \norm{x_0} \,  \omega_1^{p' -1} \,T_1                 \\
      & + 2\pi^2{(l+1)}^2L^2 \Lambda_1 \, \norm{ x_0 } \, \omega_1^{p'-1} \,T_1      \\
      & + 2 \pi^2  l^2  L \Lambda_1 \, \norm{  x_0 }  \,   \omega_1^{p'-1}           \\
      & + \pi l  L \Lambda_1 \,  \norm{  x_0  } \, \omega_1^{p'-1}  \, T_1           \\
      & + \tfrac{2}{3} \pi^2l^2L\Lambda_1 \, \norm{x_0} \, \omega_1^{p'-1}\, {T_1}   \\
      & + \pi l {L} \Lambda_1 \, \norm{x_0} \, \omega_1^{p'-1} \, T_1                \\
      & + 2 \pi^2  l^2 L \Lambda_1 \, \norm{x_0} \,  \omega_1^{p'-1}                 \\
      & +  8 \pi^2 l^2 L \Lambda_1 \, \norm{x_0}  \,  \omega_1^{p'-1} \bigr\rparen
    \end{aligned}
    \\
    \labelrel{\leq}{rel:tmp:RT1:3}
      {}
    \begin{split}
      {}&  8 (r T_1) \pi^2 {(l+1)}^2 L^2 \Lambda_1 \, \omega_1^{p'-1} \, \norm{ x_0 }
      \\
      \begin{split}
        &
        + 12 \pi^2
        {(l-1)}^2
        L^2
        \Lambda_1
        \,
        \omega_1^{p'-1}
        \,
        \norm{x_0}
        \\
        & + 6 \pi^3 {(l+1)}^2   L^2 \Lambda_1 \,  \omega_1^{p'-1} \,\norm{  x_0  }
      \end{split}
    \end{split}
    \\
    \labelrel{\leq}{rel:tmp:RT1:4}
      {} &
    \pi^2 {(l+1)}^2 L^2 \Lambda_1 \, \omega_1^{p'-1} \, \norm{ x_0 }
    \\
    \MoveEqLeft[4]
    \makebox[\displaywidth-4em][r]{
        \(
        \displaystyle
        \mathrel{\times}
        \left\lparen
        8 (t_1 -t_0)
        +12
        +6\pi
        \right\rparen
        \)
      }
  \end{align*}
  \endgroup
  where we
  \begin{steps}
    \item[\eqref{rel:tmp:RT1:1}] use the triangle inequality of norms,
    \item[\eqref{rel:tmp:RT1:2}] substitute in the individual bounds from \Cref{lemma:bounds:QV1}, \Cref{lemma:bounds:QV2}, \Cref{lemma:bounds:Q1}, \Cref{lemma:bounds:RL1}, \Cref{lemma:bounds:I}, \Cref{lemma:bounds:Q1T} and \Cref{lemma:bounds:H},
    \item[\eqref{rel:tmp:RT1:3}] overestimate contributions of \( \pi\), \(L\) and \(l\), remember the fact that \( \omega_1 \in \interval[open right]{1}{\infty}\) and  hence \( T_1 = 2\pi \omega_1^{-1} \leq 2 \pi\),
    \item[\eqref{rel:tmp:RT1:4}] remember \((rT_1) \leq t_1 - t_0 \) and sum up the expressions to conclude the statement of this lemma. \qedhere
  \end{steps}
\end{proofEnd}
\begin{theoremEnd}[%
    category=auxiliary-result,%
  ]{lemma}[Bounds on \(R\)]\label{lemma:bounds:R}
  \mbox{}
  Consider \(R\)  as defined in~\eqref{eqn:functional-expansion:R} with \Cref{ass:functional-expansion}.
  Then,
  \begin{align}
    \norm{R(t_s, t_e)}
     & \leq
    \tfrac {2 } 3 \pi^2 {(l+1)}^3 L \Lambda_1 \, \norm{x_0} \, \omega_1^{p'''-2} \, T_1.
  \end{align}
\end{theoremEnd}
\begin{proofEnd}
  The result is established by straighforwardly calculating
  \begin{align*}
    \norm{R(t_s, t_e)}
     & \labelrel{\leq}{rel:tmp:pfe:R:1}\tfrac 1 6 l {(l+1)}^2 \omega_1^{p'''} L \Lambda_1 \, \norm{x_0} \,  {(t_e - t_s)}^3
    \\
     & \labelrel{\leq}{rel:tmp:pfe:R:2}\tfrac {2 } 3 \pi^2 {(l+1)}^3 L \Lambda_1 \, \norm{x_0} \, \omega_1^{p'''-2}\,T_1,
  \end{align*}
  where we
  \begin{steps}
    \item[\eqref{rel:tmp:pfe:R:1}] use~\eqref{eqn:es-solution-bound}, \(t_s, t_e \in \interval[]{t_0'}{t_1}\),  the Lipschitz property in~\Cref{ass:fi:2:3}  and definition~\eqref{eqn:def:p:3}, then solve the integral,
    \item[\eqref{rel:tmp:pfe:R:2}] overestimate the contribution of \(l\), use the fact that \( t_e - t_s \leq T_1\) and remember \(T_1\) from ~\Cref{ass:functional-expansion} to conclude the statement of this lemma. \qedhere
  \end{steps}
\end{proofEnd}

\begin{theoremEnd}[%
    category=auxiliary-result,%
  ]{lemma}[Bounds on \( J_0\)]\label{lemma:bounds:J0}
  Consider \(J_0\)  as defined in~\eqref{eqn:proof:functional-expansion:J0}, 
  with~\Cref{ass:functional-expansion},
  and let \( t_s \in \interval[]{t_0'}{t_1}\) and  \({t_e \in \interval[]{t_s}{t_1} }\).
  Then,
  \begin{align}
    \norm{ J_0(\theta)  }
    & \leq
    L  \,  \norm{ x_{\omega_1}(\theta; t'_0, x'_0) - \bar{x}(\theta; t'_0, x'_0) }   \odif{\theta}.
  \end{align}
\end{theoremEnd}
\begin{proofEnd}
  The result is established by straighforward calculations
  using solely the uniform global Lipschitz property from~\Cref{ass:fi:2:1}.
\end{proofEnd}
\begin{theoremEnd}[%
    category=auxiliary-result,%
  ]{lemma}[Bounds on \( J_{ij}\)]\label{lemma:bounds:Jij}
  Consider \(J_{ij}\)  as defined in~\eqref{eqn:proof:functional-expansion:Jij} with \Cref{ass:functional-expansion}, and let \( t_s \in \interval[]{t_0'}{t_1}\) and  \({t_e \in \interval[]{t_s}{t_1} }\).
  Then, for all \( i \in \set{1, \ldots, l-1}\) and all \( j \in \set{ i+1,  \ldots, l}\),
  \begin{align}
    \int_{t_s}^{t_e} \norm{ J_{ij}(\theta) }  \odif{\theta}
    & \leq
    2 \pi  L \,  \norm{ x_{\omega_1}(\theta; t'_0, x'_0) - \bar{x}(\theta; t'_0, x'_0)} .
  \end{align}
\end{theoremEnd}
\begin{proofEnd}
  To analyze the expression~\( \norm{J_{ij}(\theta) } \),	three cases have to be distinguished:
  \begin{subequations}
    \begin{enumerate}
      \item If \( p_i + p_j = 1\), \( \gamma_{ij}(\omega_1)\) is a constant such that \({ |\gamma_{ij}| \leq  \pi } \).
            In that case, we estimate
            \begin{align}
              \norm{ J_{ij}(\theta) }
              &
              \leq
              2 \pi   L
              \norm{ x_{\omega_1}(\theta; t'_0, x'_0) - \bar{x}(\theta; t'_0, x'_0)}
            \end{align}
            where we
            use the uniform global Lipschitz property from~\Cref{ass:fi:2:2}.
      \item If \( p_i + p_j > 1\), this implies \( \gamma_{ij}(\omega) = 0\) or \({ \commutator{f_i}{f_j}(\theta, x) = 0}\) for all \( \theta \in \mathbb{R}\) and all \( x \in \mathbb{R}^n\)  by \Cref{ass:pi}. Accordingly,
            \begin{align}
              J_{ij}(\theta)
              & = 0.
              \label{eqn:tmp:pfe:bounds:c2}
            \end{align}
      \item If \( p_i + p_j \in \interval[open]{0}{1}\), by definition of \( \gamma_{ij}(\omega)\), we infer that \( \lim_{\omega \to \infty } \gamma_{ij}(\omega) = 0\), and hence
            \begin{align}
              \norm{ J_{ij}(\theta) }
              &
              \leq
              \pi L    \norm{ x_{\omega_1}(\theta; t'_0, x'_0) - \bar{x}(\theta; t'_0, x'_0)}
            \end{align}
            where we
            use the uniform global Lipschitz property from~\Cref{ass:fi:2:2}.
    \end{enumerate}
  \end{subequations}
  Taking the maximum value from each case completes the proof.
\end{proofEnd}

\begin{theoremEnd}[%
    category=lipschitz-proofs,%
  ]{lemma}[Lipschitz properties I]\label{lemma:lipschitz-properties-1}
  \mbox{}
  Let \( a \in \mathbb{R} \), \( {b \in \mathbb{R}\!\setminus\!\set{0}}\) and consider the functions \( f_i : \mathbb{R} \to \mathbb{R} \) given by~\eqref{eqn:example-system:0},\eqref{eqn:example-system:1} and~\eqref{eqn:example-system:2}.
  Then, for each \( i \in \{0, 1, 2\}\), the function \( f_i \) is globally uniformly {\(L_i\)-Lipschitz} with Lipschitz constants
  \begin{align*}
    L_0=\abs{a}, &  & L_1=L_2 = 3\abs{b}k^{\tfrac 1 2}
  \end{align*}
  and \(f_i(t, 0) = 0\).
\end{theoremEnd}
\begin{proofEnd}
  Function \( f_0 \) is a linear function, and the Lipschitz property holds globally with \( L_0 = \abs{a}\) for any \( a \in \mathbb{R} \).
  
  To show the claimed properties for \(f_1\) and \(f_2\), we employ the fundamental theorem of calculus.
  
  To do so, we first need to identify points, at which derivatives to \(f_1\) and \(f_2\) do not exist. The single critical point is \( x_0 = 0\).
  For all \( x \in \interval[open]{-\infty}{x_0}\), the derivatives to \( f_1\) and \( f_2\) exist and are given by
  \begin{align}
    \pdif{2} f_1(t, x) & = -bk^{\tfrac 1 2}\sin(\log(\tfrac 1 2 x^2)) - 2bk^{\tfrac 1 2} \cos(\log(\tfrac 1 2 x^2)) \\
    \pdif{2} f_2(t, x) & = -bk^{\tfrac 1 2}\cos(\log(\tfrac 1 2 x^2)) + 2bk^{\tfrac 1 2} \sin(\log(\tfrac 1 2 x^2)).
  \end{align}
  Likewise, for all  \( x \in \interval[open]{x_0}{\infty}\), the derivatives to \( f_1\) and \( f_2\) exist and are given by
  \begin{align}
    \pdif{2} f_1(t, x) & = bk^{\tfrac 1 2}\sin(\log(\tfrac 1 2 x^2)) + 2bk^{\tfrac 1 2} \cos(\log(\tfrac 1 2 x^2)) \\
    \pdif{2} f_2(t, x) & = bk^{\tfrac 1 2}\cos(\log(\tfrac 1 2 x^2)) - 2bk^{\tfrac 1 2} \sin(\log(\tfrac 1 2 x^2)).
  \end{align}
  
  Observe, that in particular \(  |\pdif{2} f_1(x)| \leq 3\abs{b}k^{\tfrac 1 2} \eqqcolon L_1 \) and \( |\pdif{2} f_2(x)| \leq 3\abs{b}k^{\tfrac 1 2} \eqqcolon L_2 \) for all \( {x \in \mathbb{R}\! \setminus\! \{x_0\}}\).
  
  Now let \( x_1, x_2 \in \mathbb{R} \), \(x_1 < x_2 \) arbitrary and calculate for \( i \in \set{1, 2}\)
  \begin{align}
    \mathrlap{|f_i(t, x_2) - f_i(t, x_1) |}  \notag \\
    & \leq
    |f_i(t, x_2) - f_i(t, x_0)| + |f_i(t, x_0) - f_i(t, x_1)|
    \label{eqn:lipschitz-properties-1:adding-zero-and-triange-inequality}
    \\
    & \leq
    \left|
    \int_{x_0}^{x_2} \pdif{2} f_i(t, x) \odif{x}
    \right|
    +
    \left|
    \int_{x_1}^{x_0} \pdif{2} f_i(t, x) \odif{x}
    \right|
    \label{eqn:lipschitz-properties-1:fundamental-theorem-of-calculus}
    \\
    & \leq
    \int_{x_0}^{x_2} |\pdif{2} f_i(t, x)| \odif{x}
    +
    \int_{x_1}^{x_0} |\pdif{2} f_i(t, x)| \odif{x}
    \label{eqn:lipschitz-properties-1:triangle-inequality}
    \\
    & \leq
    \int_{x_0}^{x_2} L_i \odif{x}
    +
    \int_{x_1}^{x_0} L_i \odif{x}
    \label{eqn:lipschitz-properties-1:derivative-bound}
    \\
    & = L_i (x_2 - x_1),
  \end{align}
  which proves the claimed global \(L_i\)-Lipschitz properties.
  Observing \(\pdif{1} f_i(t, x) = 0\) for all \(t \in \mathbb{R}\) completes the proof.
\end{proofEnd}

\begin{theoremEnd}[%
    category=lipschitz-proofs,%
  ]{lemma}[Lipschitz properties II]\label{lemma:lipschitz-properties-2}
  Consider the setup from \Cref{lemma:lipschitz-properties-1}.
  For each \( i \in \{0, 1, 2\}\) and \(j \in \{0, 1, 2\}\), the function \( \mathcal{L}_{f_j}f_i : \mathbb{R} \to \mathbb{R} \) is globally uniformly {\(L_{ij}\)-Lipschitz} with Lipschitz constants
  \begin{align}
    L_{ij} & = 6 \abs{b} k^{\tfrac 1 2} L_j + L_i L_j
    \label{eqn:L:ij}
  \end{align}
  and \(\pdif{1}\mathcal{L}_{f_j}f_i(t, 0) = 0\).
\end{theoremEnd}

\begin{proofEnd}
  Per definition,
  \begin{align}
    \mathcal{L}_{f_j}f_i(t, x) & \coloneqq \pdif{2} f_i(t, x) \cdot f_i(t, x),
  \end{align}
  that is, the function \( \mathcal{L}_{f_j}f_i \) might be undefined on critical points \( x_\text{crit}\) where either \( f_i(t, x_\text{crit}) \) or \( \pdif{2} f_j(t, x_\text{crit})\) are undefined.
  Observe, that the only such critical point is \( x_\text{crit} = x_0 = 0 \).

  For all \( x \in \interval[open]{-\infty}{x_0}\), the derivatives to \( \mathcal{L}_{f_j}f_i \) exist and are given by
  \begin{align}
    \pdif{2} \mathcal{L}_{f_j}f_i(t, x) & = \pdif{2}^2 f_i(t, x) \cdot f_j(t, x) + \pdif{2} f_i(t, x) \cdot \pdif{2} f_j(t, x).
  \end{align}
  with
  \begin{subequations}
    \begin{align}
      \pdif{2}^2 f_0(t, x) & =0 \\
      \pdif{2}^2 f_1(t, x) & = \tfrac{2}{x} bk^{\tfrac 1 2} \left(
        - \cos ( \log(\tfrac 1 2 x^2) ) +2 \sin(\log(\tfrac 1 2 x^2))
        \right) \\
      \pdif{2}^2 f_2(t, x) & = \tfrac{2}{x} bk^{\tfrac 1 2} \left(
        \sin ( \log(\tfrac 1 2 x^2) ) +2 \cos(\log(\tfrac 1 2 x^2))
        \right)
    \end{align}
  \end{subequations}
  Likewise, for all  \( x \in \interval[open]{x_0}{\infty}\), the derivatives to \(  \mathcal{L}_{f_j}f_i \) exist and are given by
  \begin{subequations}
    \begin{align}
      \pdif{2}^2 f_0(t, x) & =0 \\
      \pdif{2}^2 f_1(t, x) & = \tfrac{2}{x} bk^{\tfrac 1 2}    \left(
        \cos ( \log(\tfrac 1 2 x^2) ) -2 \sin(\log(\tfrac 1 2 x^2))
        \right) \\
      \pdif{2}^2 f_2(t,x) & = \tfrac{2}{x} bk^{\tfrac 1 2} \left(
        -\sin ( \log(\tfrac 1 2 x^2) ) -2 \cos(\log(\tfrac 1 2 x^2))
        \right)
    \end{align}
  \end{subequations}
  Observe, that in particular \(  |\pdif{2}^2 f_i(t, x)| \leq \frac{6}{|x|}\abs{b} k^{\tfrac 1 2} \) for all \({ i \in \set{0, 1, 2}}\) and all \( { x \in \mathbb{R} \setminus \{x_0\} }\).
  Then, for all \({ i,j \in \set{0, 1, 2}}\), the estimate
  \begin{align}
    |\pdif{2} \mathcal{L}_{f_j}f_i(t, x)| & \leq 6\abs{b} k^{\tfrac 1 2} L_j + L_i L_j
  \end{align}
  is satisfied for all \( x \in \mathbb{R} \!\setminus\! \{x_0\}\).
  Now let \( x_1, x_2 \in \mathbb{R} \), \(x_1 < x_2 \) arbitrary and calculate for \( i,j \in \{0, 1, 2\}\)
  \begin{subequations}
    \begin{align}
      \mathrlap{\abs{\mathcal{L}_{f_j}f_i(t, x_2) - \mathcal{L}_{f_j}f_i(t, x_1) }}
      \phantom{aaaaaaaaaa}
      \notag\\
      \begin{split}
        &{} \leq  \mathrel{\phantom{+}} \abs{\mathcal{L}_{f_j}f_i(t, x_2) - \mathcal{L}_{f_j}f_i(t, x_0)} \\
        &{} \mathrel{\phantom{\leq}}
        +
        \abs{\mathcal{L}_{f_j}f_i(t, x_0) - \mathcal{L}_{f_j}f_i(t, x_1)}
      \end{split}
      \label{eqn:lipschitz-properties-2:adding-zero-and-triange-inequality}
      \\
      \begin{split}
        &{} \leq \mathrel{\phantom{+}}
        \abs{
            \int_{x_0}^{x_2} \pdif{2} \mathcal{L}_{f_j}f_i(t, x) \odif{x}
          }
        \\
        &{} \mathrel{\phantom{\leq}}
        +
        \abs{\int_{x_1}^{x_0} \pdif{2}\mathcal{L}_{f_j}f_i(t, x) \odif{x}
          }
      \end{split}
      \label{eqn:lipschitz-properties-2:fundamental-theorem-of-calculus}
      \\
      \begin{split}
        &{} \leq \mathrel{\phantom{+}}
        \int_{x_0}^{x_2} \abs{\pdif{2}\mathcal{L}_{f_j}f_i(t, x)} \odif{x}
        \\
        &{} \mathrel{\phantom{\leq}}
        +
        \int_{x_1}^{x_0} \abs{\pdif{2}\mathcal{L}_{f_j}f_i(t, x)} \odif{x}
      \end{split}
      \label{eqn:lipschitz-properties-2:triangle-inequality}
      \\
      \begin{split}
        &{} \leq \mathrel{\phantom{+}}
        \int_{x_0}^{x_2} 6\abs{b} k^{\tfrac 1 2}L_j +L_i L_j \odif{x}
        \\
        &{} \mathrel{\phantom{\leq}}
        +
        \int_{x_1}^{x_0} 6 \abs{b}k^{\tfrac 1 2}L_j +L_i L_j \odif{x}
      \end{split}
      \label{eqn:lipschitz-properties-2:final}
      \\
      & = (6\abs{b} k^{\tfrac 1 2} L_j+ L_i L_j) (x_2 - x_1),
    \end{align}
    which with the definition~\eqref{eqn:L:ij}
  \end{subequations}
  proves the claimed global {\( L_{ij} \)-Lipschitz} properties.
  Observing that \(\pdif{1}\mathcal{L}_{f_j}f_i(t, x) = 0\) for all \(t\) and for all \(x \in \mathbb{R}\) completes the proof.
\end{proofEnd}

\begin{theoremEnd}[%
    category=lipschitz-proofs,%
  ]{lemma}[Lipschitz properties III]\label{lemma:lipschitz-properties-3}
  Consider the setup from \Cref{lemma:lipschitz-properties-1}.
  For each \( i, j, m \in \set{0, 1, 2}\), the function \( \mathcal{L}_{f_m}\mathcal{L}_{f_j}f_i : \mathbb{R} \to \mathbb{R} \) is globally uniformly {\(L_{ijm}\)-Lipschitz} with Lipschitz constants
  \begin{align}
    L_{ijm} & = 16 \abs{b}k^{\tfrac 1 2} L_j L_m + 6 \abs{b}k^{\tfrac 1 2} L_i L_m + L_i L_j L_m.
    \label{eqn:L:ijm}
  \end{align}
\end{theoremEnd}
\begin{proofEnd}
  Per definition,
  \begin{align}
    \mathcal{L}_{f_m}\mathcal{L}_{f_j}f_i(t, x) & \coloneqq \pdif{2}(\mathcal{L}_{f_j}f_i(t, x)) \cdot f_m(t, x),
  \end{align}
  that is, the function \(\mathcal{L}_{f_m}\mathcal{L}_{f_j}f_i \) might be undefined on critical points \( x_\text{crit}\) where either \( \pdif{2} (\mathcal{L}_{f_j}f_i(t, x_\text{crit})) \) or \( f_m(t, x_\text{crit})\) is undefined.
  Observe, that the single such critical point is \( x_\text{crit} = x_0 = 0 \).
  
  For all \( x \in \interval[open]{-\infty}{x_0}\), the partial derivatives to \( \mathcal{L}_{f_m}\mathcal{L}_{f_j}f_i \) exist and are given by\\
  \vspace*{-4ex}
  \begin{subequations}
    \begin{align}
      \pdif{2}  \mathcal{L}_{f_m}\mathcal{L}_{f_j}f_i(t, x)
      \begin{split}
         & {} = \phantom{+} \pdif{2}^2 \mathcal{L}_{f_j}f_i(t, x) \cdot f_m(t, x)
        \\
         & {} \mathrel{\phantom{=}}	+
        \pdif{2} \mathcal{L}_{f_j}f_i(t, x) \cdot \pdif{2} f_m(t, x)
      \end{split}
      \\
      \begin{split}
         & {} =  \phantom{+}\pdif{2}^3 f_i(t, x) \cdot f_j(t, x) \cdot f_m(t, x) \\
         & {} \mathrel{\phantom{=}}                                              %
        +
        \pdif{2}^2 f_i(t, x) \cdot \pdif{2} f_i(t, x) \cdot f_m(t, x)
        \\
         & {} \mathrel{\phantom{=}} +
        \pdif{2}^2 f_i(t, x) \cdot f_j(t, x) \cdot \pdif{2} f_m(t, x)
        \\
         & {} \mathrel{\phantom{=}} +
        \pdif{2} f_i(t, x) \cdot \pdif{2} f_j(t, x) \cdot \pdif{2} f_m(t, x)
      \end{split}
      \label{eqn:lipschitz-properties-3:partial-derivative-mji}
    \end{align}
  \end{subequations}
  with \( \pdif{2} f_i(t, x) \) from the proof of \Cref{lemma:lipschitz-properties-1}, \( \pdif{2}^2 f_i(t, x)\) from the proof of \Cref{lemma:lipschitz-properties-2} and
  \vspace{-1ex}
  \begin{subequations}
    \begin{align}
      \pdif{2}^3 f_0(x) & =0                                                                                          \\
      \pdif{2}^3 f_1(x) & = \tfrac{2}{x^2}b k^{\tfrac 1 2} \left( \phantom{-} 5 \cos ( \log(\tfrac 1 2 x^2) ) \right) \\
      \pdif{2}^3 f_2(x) & = \tfrac{2}{x^2}b k^{\tfrac 1 2} \left( -5\sin ( \log(\tfrac 1 2 x^2) ) \right)
    \end{align}
  \end{subequations}
  for all \( i \in \set{0,1,2}\).
  Likewise, for all  \( x \in \interval[open]{x_0}{\infty}\), the derivatives to \(   \mathcal{L}_{f_m}\mathcal{L}_{f_j}f_i  \) exist and are given by~\eqref{eqn:lipschitz-properties-3:partial-derivative-mji} with \( \pdif{2} f_i(x) \) from the proof of \Cref{lemma:lipschitz-properties-1}, \( \pdif{2}^2 f_i(x)\) from the proof of \Cref{lemma:lipschitz-properties-2} and
  \vspace{-3ex}
  \begin{subequations}
    \begin{align}
      \pdif{2}^3 f_0(x)   & =0                                                                                       \\
      \pdif{2}^3 f_1(x)   & = \tfrac{2}{x^2}b k^{\tfrac 1 2} \left( -5\cos ( \log(\tfrac 1 2 x^2) )\right)           \\
      \pdif{2}^3 f_2(x^2) & = \tfrac{2}{x^2}b k^{\tfrac 1 2} \left(\phantom{-}5\sin ( \log(\tfrac 1 2 x^2) ) \right)
    \end{align}
  \end{subequations}
  for all \( i \in \set{0,1,2}\).
  Observe, that in particular \({ |\pdif{2}^3 f_i(x)| \leq \tfrac{10}{{|x|}^2} \abs{b} k^{\tfrac 1 2} }\)
  for all \( { i \in \set{0, 1, 2} }\) and all \( { x \in \mathbb{R} \!\setminus\! \set{x_0} }\).
  Consequently, for all \( i,j,m \in \set{0, 1, 2}\), the estimate
  \begin{align}
    \begin{split}
      \mathrlap{\abs{\pdif{2} \mathcal{L}_{f_m}\mathcal{L}_{f_j}f_i(x)}}
      \phantom{aaa}
      \notag                                                                     \\
       & {} \leq
      \frac{10}{{|x|^2}} \abs{b}k^{\tfrac 1 2} \cdot L_j \, |x| \cdot L_m |x|
      \\
       & {} \phantom{=} + \frac{6}{|x|}\abs{b}k^{\tfrac 1 2} L_i \cdot L_m |x|
      \\
       & {} \phantom{=} + \frac{6}{|x|}\abs{b}k^{\tfrac 1 2} L_j |x| \cdot  L_m
      + L_i L_j L_m
    \end{split}
    \\
     & = 16\abs{b}k^{\tfrac 1 2} L_j L_m + 6\abs{b}k^{\tfrac 1 2} L_i L_m + L_i L_j L_m
  \end{align}
  is satisfied for all \( x \in \mathbb{R} \! \setminus \! \set{x_0}\).
  Now let \( x_1, x_2 \in \mathbb{R} \), \(x_1 < x_2 \) arbitrary and calculate for \( i,j \in \{0, 1, 2\}\)
  \begin{align}
    \begin{split}
      \mathrlap{\abs*{\mathcal{L}_{f_m}\mathcal{L}_{f_j}f_i(x)(t, x_2) - \mathcal{L}_{f_m}\mathcal{L}_{f_j}f_i(x)(t, x_1) }}
      \notag
      \\
       & \leq
      \abs*{\mathcal{L}_{f_m}\mathcal{L}_{f_j}f_i(x)(t, x_2) - \mathcal{L}_{f_m}\mathcal{L}_{f_j}f_i(t, x_0)}                           \\
       & \mathrel{\phantom{=}}{} + \abs*{\mathcal{L}_{f_m}\mathcal{L}_{f_j}f_i(t, x_0) - \mathcal{L}_{f_m}\mathcal{L}_{f_j}f_i(t, x_1)}
    \end{split}
    \\
    \begin{split}
       & \leq
      \abs{
          \int_{x_0}^{x_2} \pdif{2}\mathcal{L}_{f_m}\mathcal{L}_{f_j}f_i(t,x) \odif{x}
        }
      \\
       & \mathrel{\phantom{=}} {} +
      \abs{
          \int_{x_1}^{x_0} \pdif{2}\mathcal{L}_{f_m}\mathcal{L}_{f_j}f_i(t, x) \odif{x}
        }
    \end{split}\label{eqn:lipschitz-properties-3:fundamental-theorem-of-calculus}
    \\
    \begin{split}
       & \leq
      \int_{x_0}^{x_2} \abs{\pdif{1}\mathcal{L}_{f_m}\mathcal{L}_{f_j}f_i(x)(x)} \odif{x}
      \\
       & \mathrel{\phantom{=}} {}+
      \int_{x_1}^{x_0} \abs{\pdif{1}\mathcal{L}_{f_m}\mathcal{L}_{f_j}f_i(x)} \odif{x}
    \end{split}
    \label{eqn:lipschitz-properties-3:triangle-inequality}
    \\
    \begin{split}
       & \leq
      \int_{x_0}^{x_2} 16 \abs{b}k^{\tfrac 1 2} L_j L_m + 6 \abs{b}k^{\tfrac 1 2} L_i L_m + L_i L_j L_m \odif{x}
      \\
       & \mathrel{\phantom{=}} {} +
      \int_{x_1}^{x_0} 16 \abs{b}k^{\tfrac 1 2}  L_j L_m + 6 \abs{b}k^{\tfrac 1 2} L_i L_m + L_i L_j L_m \odif{x}
    \end{split}
    \label{eqn:lipschitz-properties-3:partial-derivative-bounds}
    \\
     & = (16 \abs{b}k^{\tfrac 1 2}L_j L_m + 6 \abs{b}k^{\tfrac 1 2}L_i L_m + L_i L_j L_m) (x_2 - x_1),
  \end{align}
  which with the definition~\eqref{eqn:L:ijm} proves the claimed global {\( L_{ijm} \)-Lipschitz} properties.
\end{proofEnd}

\section{Proofs}\label{sec:proofs}
 \subsection{Proofs of Main Results}
  
\begin{proof}[of \Cref{lemma:approximation-lemma}]\label{proof:approximation-lemma}
  In virtue of \Cref{ass:ui,ass:fi},
  the following are satisfied:
  \begin{enumerate}
    \item The pair \((\mathbb{R}^n, \norm{\cdot})\) forms a Banach space.
    \item The vector field			     \( f_{\omega}(t, x)  \coloneqq  \sum_{i=0}^l \omega^{p_i} f_i(t, x) \cdot u_i(\omega t) \)
          is continuous, and uniformly globally Lipschitz on \( \mathbb{R} \times \mathbb{R}^n\).
  \end{enumerate}
  Consequently, the Picard-Lindel{ö}f Theorem guarantees, that for every \( t_0 \in \mathbb{R} \), \( x_0 \in \mathbb{R}^n \) and \( \omega \in \interval[open]{0}{\infty}\), there exists a \( t_u(t_0, x_0, \omega) \in \interval[open]{0}{\infty}  \) such that the trajectory of system~\eqref{eqn:es-system} with initial condition \( x(t_0) = x_0 \) exists uniquely on \( \interval[]{t_0}{t_0 + t_u(t_0, x_0, \omega)} \).
  Let \( t_m(t_0, x_0, \omega)\) be the supremal value of \( t_u(t_0, x_0, \omega)\).
  Since \(f_\omega\) is uniformly \emph{globally} Lipschitz, the solution exists globally and \( {t_m = \infty}\).

  Observe, that both solutions \( x_{\omega}\) and \( \bar{x}\) exist {(at least)} on their common interval \(  \mathcal{T}(t_0, x_0, \omega) \coloneqq \interval[]{t_0}{t_0+t_f} \cap \interval[]{t_0}{t_0+ t_m(t_0, x_0, \omega)} = \interval[]{t_0}{t_0+t_f}\) per assumptions of this theorem.

  To prove the result, let us proceed by contradiction. Namely, let us falsify the following statement obtained by negating a portion of the statement of this theorem:

  \begin{customassumption}[Contradiction Statement]\label{ass:contradiction-assumption}
    Let \( D \in \interval[open]{0}{\infty}\) be arbitrary but fixed.
    Then for each \( \omega^\star \in \interval[open]{0}{\infty}\), there exists an \({ \omega_1 \in \interval[open]{\omega^\star}{\infty}}\),  \({ t'_0 \in \mathbb{R}} \), a~\( {x'_0 \in {\mathbb{R}^n \! \setminus \!\set{0}}} \) and  \({ t_s \in \mathcal{T}(t'_0, x'_0, \omega_1)}\), such that the trajectories of systems~\eqref{eqn:es-system} and~\eqref{eqn:lbs} through \({ x(t'_0) = x'_0 }\) and \( {\bar{x}(t'_0) = x'_0}\) satisfy
    \par
    \vspace{\abovedisplayskip}
    \hfill \(\displaystyle \norm{ x_{\omega_1}(t_s; t'_0, x'_0) - \bar{x}(t_s; t'_0, x'_0) } \geq D \,\norm{x_0'} \) \;. \qedhere
  \end{customassumption}

  Let \( D \) be arbitrary but fixed. For (\nameref{ass:contradiction-assumption}) to hold for all \( \omega^\star \in \interval[open]{0}{\infty}\), it must hold in particular for \( \omega^\star\) as defined in~\eqref{eqn:def:omega-star}.

  Note that \( x_{\omega_1}(t; t'_0, x'_0) \) and \( \bar{x}(t; t'_0, x'_0) \) are absolutely continuous on \( \interval[]{t'_0}{t_s} \).
  According to (\nameref{ass:contradiction-assumption}), having \( \norm{ x_{\omega_1}(t_s; t'_0, x'_0) - \bar{x}(t_s; t'_0, x'_0) } \geq D \norm{x_0'} \) with \({ \norm{ x_{\omega_1}(t'_0; t'_0, x'_0) - \bar{x}(t'_0; t'_0, x'_0) }  = 0} \) implies the existence of a \( t_2 \in \interval[open left]{t_0'}{t_s} \), such that
  \begin{subequations}
    \begin{align}
      \norm{ x_{\omega_1}(t_2; t'_0, x'_0) - \bar{x}(t_2; t'_0, x'_0) } & = D \, \norm{x_0'} 
      \label{eqn:tmp:pfe:ca:1}
    \end{align}
    and
    \begin{align}
      \norm{ x_{\omega_1}(t; t'_0, x'_0) - \bar{x}(t; t'_0, x'_0) } & < D \, \norm{x_0'}
      \label{eqn:tmp:pfe:ca:2}
    \end{align}
    for all \( t \in \interval[open right]{t'_0}{t_2} \).
  \end{subequations}
  To reach a contradiction, we will now examine the distance between \(x_{\omega_1}(t; t'_0, x'_0)\) and \(\bar{x}(t; t'_0, x'_0)\).
  Let \( t_1 \in \interval[]{t'_0}{t_2} \) be arbitrary but fixed.
  Since \( \bar{x}(t; t'_0, x'_0)\) exists on \( \interval[]{t'_0}{t_1}\), we may express the solution to~\eqref{eqn:lbs} as
  \begin{align}
    \begin{split}
      \MoveEqLeft
      \bar{x}(t_1; t'_0, x'_0)
      ={}  x'_0 + \int_{t'_0}^{t_1} f_0(\theta, \bar{x}(\theta; t'_0, x'_0))
      \\
       & + \sum_{\mathclap{\substack{1\leq i <l \\i<j \leq l}}} \;
      \lim_{\omega \to \infty}
      \commutator{f_i}{f_j}(\theta, \bar{x}(\theta; t'_0, x'_0)) \gamma_{ij}(\omega) \odif{\theta}.
    \end{split}
    \label{eqn:A68}
  \end{align}
  Furthermore, since the conditions of \Cref{lemma:functional-expansion} are met with the selected values of \(\omega_1\), \(t_0 = t'_0\), \(t_1\) and \(x_0 = x'_0\) and the choice of \(  D \), we may express the solution to~\eqref{eqn:es-system} as
  \begin{align}
    \begin{split}
      \MoveEqLeft
      x_{\omega_1}(t_1; t'_0, x'_0)
      =
      x'_0
      + \int_{t'_0}^{t_1} f_0(\theta, x_{\omega_1}(\theta; t'_0, x'_0))
      \\
       & + \sum_{\mathclap{\begin{limitarray}
                                     1 &\leq i &< l \\
                                     i&<j&\leq l
                                   \end{limitarray}}} \commutator{f_i}{f_j}(\theta, x_{\omega_1}(\theta; t'_0, x'_0)) \gamma_{ij}(\omega) \odif{\theta} \\
       & + \sum_{q=0}^{r-1} \!R(t'_q, t'_{q+1}) + \!R(t'_b, t_1) + \!R_{T_1}(t'_0, t_1),
    \end{split}
    \label{eqn:A69}
  \end{align}
  with \( T_1\), \( r  \), \( t_q\), \(R\) and \( R_{T_1}\) defined in \Cref{ass:functional-expansion}. Combining~\eqref{eqn:A68} with~\eqref{eqn:A69}, we obtain with
  \begin{subequations}
    \begin{align}
      J_0(\theta)     \coloneqq{} & f_0(\theta, x_{\omega_1}(\theta; t'_0, x'_0)) - f_0(\theta, \bar{x}(\theta; t'_0, x'_0))
      \label{eqn:proof:functional-expansion:J0}
      \\
      J_{ij}(\theta)               \coloneqq {}
      \begin{split}
         & \commutator{f_i}{f_j}(\theta, x_{\omega_1}(\theta; t'_0, x'_0)) \,  \gamma_{ij}(\omega_1)
        \\
         &
        - \lim_{\omega \to \infty} \commutator{f_i}{f_j}(\theta, \bar{x}(\theta; t'_0, x'_0)) \, \gamma_{ij}(\omega)
      \end{split}
      \label{eqn:proof:functional-expansion:Jij}
    \end{align}
  \end{subequations}
  the expression
  \begingroup
  \allowdisplaybreaks
  \begin{align}
    \mathrlap{\norm{ x_{\omega_1} - \bar{x}}(t_1; t'_0, x'_0)}
    \phantom{\leq{}}
    \notag
    \\
    \labelrel{\leq}{rel:tmp:e:1}
      {}
    \begin{split}
       &
      \mathrel{\phantom{+}}
      \int_{t'_0}^{t_1}
      \norm{ J_0(\theta) } \odif{\theta}
      +
      \int_{t'_0}^{t_1}
      \sum_{\mathclap{\begin{limitarray}
                                1 &\leq i &< l \\
                                i&<j&\leq l
                              \end{limitarray}}}
      \norm{J_{ij}(\theta) }
      \odif{\theta}                                    \\[-10pt]
       & + \sum_{q=0}^{r-1} \!\norm{R(t'_q, t'_{q+1})}
      \!+ \! \norm{R(t'_r, t_1)}
      \!+ \! \norm{R_{T_1}(t'_0, t_1) }
      \label{eqn:proof:new-approximation-lemma:23}
    \end{split}
    \\
    \labelrel{\leq}{rel:tmp:e:2}
      {}
    \begin{split}
       &
      \mathrel{\phantom{+}}
      L  \int_{t_0'}^{t_1}  \norm{ x_{\omega_1}(\theta; t'_0, x'_0) - \bar{x}(\theta; t'_0, x'_0) }  \odif{\theta}
      \\[-6pt]
       &
      +
      l(l-1) \pi  L \int_{t_s}^{t_e}  \norm{ x_{\omega_1}(\theta; t'_0, x'_0) - \bar{x}(\theta; t'_0, x'_0)} \odif{\theta} \\
       & +  \pi^2 {(l+1)}^3 L^2 D  \omega_1^{p'''-2} \norm{x_0'} (t_1 - t_0')                                              \\
       &
      +\pi^2 {(l+1)}^2 L^2 D \omega_1^{p'-1} \norm{ x_0' }
      (
      8 (t_1 -t_0')
      +12
      +6\pi
      )                                                                                                                    \\
       & + 2 \pi^3 {(l+1)}^3 L^2 D  \omega_1^{p'''-2} \norm{x_0'}
    \end{split}
    \\
    \labelrel{\leq}{rel:tmp:e:3}
      {}
    \begin{split}
       &
      \mathrel{\phantom{+}}
      2 \pi l^2 L \int_{t_0'}^{t_1} \norm{ x_{\omega_1}(\theta; t'_0, x'_0) - \bar{x}(\theta; t'_0, x'_0) } \odif{\theta}
      \\
       &
      +\pi^2 {(l+1)}^3 \! L^2 \!D \omega_1^{- p^\star} \!\! \norm{ x_0' }
      (
      9 (t_1 -t_0')
      \!+\!12
      \!+8\pi
      )
    \end{split}
    \label{eqn:proof:trajectory-approximation:error-bound}
  \end{align}
  \endgroup
  describing the norm of the error between solutions, where we
  \begin{steps}[]
    \item[\eqref{rel:tmp:e:1}] use the triangle inequality to bound expressions separately,
    \item[\eqref{rel:tmp:e:2}] use the bounds from \Cref{lemma:bounds:RT1,lemma:bounds:R,lemma:bounds:J0,lemma:bounds:Jij} respectively, and finally
    \item[\eqref{rel:tmp:e:3}] sum the terms up and use the definition of \( p^\star\)  in~\eqref{eqn:def:p-star}.
  \end{steps}
  Since~\eqref{eqn:proof:trajectory-approximation:error-bound} holds for every \( t_1 \in \interval[]{t'_0}{t_2}\), we may apply Grönwall's lemma to obtain
  \begin{align}
    \mathrlap{\norm{x_{\omega_1}(t_1; t'_0, x'_0) - \bar{x}(t_1; t'_0, x'_0)}}
    \phantom{aaaaaa}{} \notag
    \\
     &                                                                       
    \leq
    \pi^2 {(l+1)}^3 L^2 D \omega_1^{- p^\star} \norm{ x_0' } \, \cdot \notag \\
     & \phantom{\leq{} +}
    (
    9 (t_1 -t_0')
    +12
    +8\pi
    )
    e^{2 \pi l^2 L (t_1 - t'_0) }
  \end{align}
  for all \( t_1 \in \interval[]{t'_0}{t_2}\).
  Considering \( t_1 = t_2\) in particular, and remembering that \( t_2 - t'_0 \in \interval[open left]{0}{t_f}\) and \( \omega_1 \in \interval[open]{\omega^\star}{\infty}\), with \( \omega^\star \) defined previously, we have that
  \begin{align}
    \norm{x_{\omega_1}(t_2; t'_0, x'_0) - \bar{x}(t_2; t'_0, x'_0)} & < D \, \norm{x'_0},
  \end{align}
  which contradicts the definition of \( t_2\) in~\eqref{eqn:tmp:pfe:ca:2}. Hence, (\nameref{ass:contradiction-assumption}) is a false statement, concluding the result.
\end{proof}

\begin{proof}[of \Cref{lemma:global-exponential-stability-result}]\label{proof:global-exponential-stability-result}
  In this proof, we take the approximation property on finite intervals from \Cref{lemma:approximation-lemma}  and extend it to the infinite interval.


  With \({ \alpha \in \interval[open right]{1}{\infty}} \) and~\( \beta \in \interval[open]{0}{\infty} \) from the assumptions of this theorem,  let \( t_f \in \interval[open]{0}{\infty} \) and \(D \in \interval[open]{0}{\infty} \) be fixed and satisfying
  \begin{align}
    \bar \alpha e^{-\bar\beta t_f} +D \in \interval[open]{0}{1}.
    \label{eqn:fb1:choice-of-tf-and-d}
  \end{align}
  We observe that any time instant \( t \in \interval[open right]{t_0}{\infty} \) may be expressed as \( t = \tau + n t_f + t_0 \) with
  \vspace*{-0.5\multicolsep}
  \begin{multicols}{2}
    \noindent
    \begin{equation}
      n \coloneqq \floor*{\tfrac{t-t_0}{t_f}}
    \end{equation}
    \begin{equation}
      \tau  \coloneqq t - t_0 - n t_f
      \vphantom{\floor*{\tfrac{t-t_0}{t_f}}}
      \label{eqn:proof:lemma:ges:tau}
    \end{equation}
  \end{multicols}
  \vspace*{-0.75\multicolsep}
  \noindent
  and hence \( \tau \in \interval[open right]{0}{t_f} \).
  Then, for every \( t \in \interval[open right]{t_0}{\infty} \), there exists a unique \( t_n \coloneqq t_0 + n t_f \).
  Since \(   n  \leq \frac{t-t_0}{t_f} \leq (n+1) \) by construction,  for every \( \delta' \in \interval[open]{0}{\infty} \), we have
  \begin{align}
    e^{-\delta' (n+1)} & \leq e^{-\delta' \tfrac{t-t_0}{t_f}} \leq e^{-\delta' n}.
    \label{eqn:fb1:star}
  \end{align}
  In view of the assumptions of this theorem,~\Cref{lemma:approximation-lemma} holds in particular for initial time \( t_n \in \interval[open right]{t_0}{\infty} \) and initial condition~\( {x_n = x_{\omega}(t_n; t_0, x_0)}\) with \(\bar\Lambda = \bar\alpha\), and the chosen pair of time horizon~\(t_f\) and approximation scale~\(D\) from~\ref{eqn:fb1:choice-of-tf-and-d}, yielding an \(\omega^\star \in \interval[open]{0}{\infty}\), such that the approximation property~\eqref{eqn:lemma:approximation-lemma:approximation} holds for all \(\omega \in \interval[open]{\omega^\star}{\infty}\) and all \( t \in \interval{t_n}{t_{n+1}}\).

  As an intermediate step, we prove by induction that under the conditions of this theorem,
  \begin{description}[
      labelindent=0pt,
      font=\normalfont\footnotesize,
      leftmargin=0pt,
      topsep=0pt,
      partopsep = 0pt,
    ]
    \item[(Induction hypothesis)] \begin{align}
            \norm{x_n} & \leq \norm{x_0} \, {(\bar \alpha e^{-\bar\beta t_f} + D )}^{n}. 
            \label{eqn:proof:lemma:ges:ih}
          \end{align}
    \item[(Base case)]Let \( n = 0 \), then trivially \( \norm{x_0} \leq  \norm{x_0} \cdot 1 \) holds.
    \item[(Induction step)] Starting with \( x_{n+1} \), we see that
          \begin{align*}
            \norm{x_{n+1}}
             & \labelrel{\leq}{rel:tmp:ges:I:1} \bar \alpha \, \norm{x_n} \, e^{-\bar\beta t_f} + D \, \norm{x_n} = \norm{x_n} \, {(\bar \alpha e^{-\bar\beta t_f }+ D)} \\
             & \labelrel{\leq}{rel:tmp:ges:I:2} \norm{x_0} \, {{(\bar \alpha e^{-\bar\beta t_f} + D)}^{n+1}},
          \end{align*}
          where we
          \begin{steps}
            \item[\eqref{rel:tmp:ges:I:1}] use the exponential stability of~\eqref{eqn:lbs} and the approximation property on finite intervals~\eqref{eqn:lemma:approximation-lemma:approximation} from \Cref{lemma:approximation-lemma},
            \item[\eqref{rel:tmp:ges:I:2}] then insert the induction hypothesis~\eqref{eqn:proof:lemma:ges:ih},
          \end{steps}
          which yields the claimed hypothesis and concludes the induction step.
  \end{description}

  Finally, we are ready to prove the statement of this theorem. A sequence of calculations shows
  \begin{subequations}
    \begin{align*}
      \norm{x(t; t_0, x_0)}
       & \labelrel{\leq}{rel:tmp:lemma:ges:1} \norm{x_n} \, {( \bar \alpha  e^{-\bar\beta (t-t_n)} +D)}                                                               \\
       & \labelrel{=}{rel:tmp:lemma:ges:2} \norm{x_n} {( \bar \alpha e^{-\bar\beta \tau} +D)}                                                                         \\
       & \labelrel{\leq}{rel:tmp:lemma:ges:3} \norm{x_0} \,  {(\bar \alpha e^{-\bar\beta t_f} +D)}^{n} {(\bar \alpha e^{-\bar\beta \tau} + D)}                        \\
       & =  \norm{x_0} \,   {(\bar \alpha e^{-\bar\beta t_f} + D)}^{n+1}  \smash[b]{\tfrac{(\bar \alpha e^{-\bar\beta \tau} + D)}{\bar \alpha e^{-\bar\beta t_f}+ D}} \\
       & \labelrel{\leq}{rel:tmp:lemma:ges:5} \alpha \, \norm{x_0} \,  e^{-\beta (t-t_0)}
    \end{align*}
  \end{subequations}
  the claimed global uniform exponential stability with \( \alpha \) and \( \beta \) defined as
  \begin{align}
    \alpha & \coloneqq \frac{ \bar \alpha +D}{\bar \alpha e^{-\bar\beta t_f}+D} \qquad\text{and}\label{eqn:lemma:ges:alpha} \\
    \beta  & \coloneqq -\frac{1}{t_f} \ln(\bar \alpha e^{-\bar\beta t_f}+D), \label{eqn:lemma:ges:beta}
  \end{align}
  respectively, where we use
  \begin{steps}
    \item[\eqref{rel:tmp:lemma:ges:1}] the exponential stability of~\eqref{eqn:lbs} and the approximation property on finite intervals from~\Cref{lemma:approximation-lemma},
    \item[\eqref{rel:tmp:lemma:ges:2}] definitions of \( t \), \( t_n\) and \( \tau \) from~\eqref{eqn:proof:lemma:ges:tau},
    \item[\eqref{rel:tmp:lemma:ges:3}]  inequality~\eqref{eqn:proof:lemma:ges:ih} from the previous inductive scheme and
    \item[\eqref{rel:tmp:lemma:ges:5}] the fact that \(e^{-\bar \beta \tau} \leq 1\),~\eqref{eqn:fb1:star} and definitions~(\ref{eqn:lemma:ges:alpha}-\ref{eqn:lemma:ges:beta}).\qedhere
  \end{steps}
\end{proof}

\begin{proof}[of \Cref{lemma:frequency-adaptation}]\label{proof:frequency-adaptation}
  The Lie-bracket system of form~\eqref{eqn:lbs} associated with~\eqref{eqn:augmented-system:state} remains the same for all \( w_k \in \interval[open]{0}{\infty}\).
  Since the origin is globally exponentially stable with respect to~\eqref{eqn:lbs} associated with~\eqref{eqn:augmented-system:state} due to the assumptions of this lemma, in view of \Cref{lemma:global-exponential-stability-result}, there exists a \(\omega^\star \in \interval[open]{0}{\infty}\), such that for all \( w_k \in \interval[open]{\omega^\star}{\infty}\), the origin is globally uniformly  exponentially stable with respect to~\eqref{eqn:augmented-system:state}.

  Consider a solution \(w(t; t_0, w_0)\) to~\eqref{eqn:augmented-system:frequency}, which is monotonically increasing by construction. Hence, either \(w(t; t_0, w_0)\) grows unbounded, or \(w(t; t_0, w_0)\) is bounded.

  Assume with contradiction intent, that \(w(t; t_0, w_0)\) grows unbounded.
  As solutions to~\eqref{eqn:augmented-system:frequency} exist globally due to the global Lipschitz properties from \Cref{ass:fi}, there exists a \( k^\star \in \mathbb{N}\), such that \( w(t; t_k, w_k) \geq \omega^\star\) for all \( k \geq k^\star\).
  Hence, the solutions \( x_{w_{k}}(t; t_{k}, x_{k})\) exponentially converge to zero for \( k \geq k^\star\), and \( x_{k^\star}\) is finite.
  Consequently, it can be shown by simple integration of~\eqref{eqn:augmented-system:frequency}, that \(w(t; t_{k^\star}, x_{k^\star})\) has a finite limit, which contradicts our initial assumption.

  Hence, finite escape time of \(w(t; t_0, w_0)\) is excluded and the statement of this lemma follows.
\end{proof}


  
\subsection{Proofs of Auxiliary Results}
  \printProofs[auxiliary-result]

\subsection{Proofs of Example Results}
  \printProofs[lipschitz-proofs]

\end{document}